\documentclass[12pt,english]{amsart}

\usepackage{amssymb, amsmath}
\usepackage{float,epsfig}
\usepackage{epsfig,amsfonts,amsbsy,amssymb,amsthm,here}
\usepackage[usenames,dvipsnames]{color}
\usepackage{graphics}
\usepackage{graphicx}
\usepackage{epsfig}
\usepackage{url}
\usepackage{arydshln,leftidx,mathtools}
\usepackage{rotating}
\usepackage{lscape}
\usepackage{setspace}
\usepackage{algpseudocode}
\usepackage{algorithm}
\usepackage[all]{xy}

\usepackage[paperwidth=220mm,paperheight=297mm,inner=3cm,outer=3cm,top=2.5cm,bottom=2.6cm]{geometry}

\newtheorem{theorem}{\bf Theorem}[section]
\newtheorem{proposition}[theorem]{\bf Proposition}

\newtheorem{corollary}[theorem]{\bf Corollary}

\newtheorem{example}[theorem]{\bf Example}

\newtheorem{remark}[theorem]{\bf Remark}
\newtheorem{lemma}[theorem]{\bf Lemma}

\newcommand{\ba}{\begin{array}}
\newcommand{\ea}{\end{array}}
\newcommand{\beano}{\begin{eqnarray*}}
\newcommand{\eeano}{\end{eqnarray*}}

\def \bmatrix#1{\left[ \begin{matrix} #1 \end{matrix} \right]}

\newcommand \noin{\noindent}

\newcommand \R{{\mathbb R}}
\newcommand \C{{\mathbb C}}
\newcommand \K{{\mathbb K}}
\newcommand \E{{\mathbf E}}

\def \P{{\mathbb P}}

\def \H{{\mathbb H}}

\newcommand \B{{\mathcal B}}
\newcommand \N{{\mathbb N}}

\newcommand \diag{\mathrm{Diag}}
\newcommand \eig{\mathrm{eig}}
\newcommand \pf{{\bf Proof: }}

\newcommand \x{{\mathbf{x}}}
\newcommand \y{{\mathbf{y}}}
\newcommand \z{{\mathbf{z}}}
\newcommand \0{{\mathbf{0}}}
\def \a{{\mathbf{a}}}
\def \b{{\mathbf{b}}}

\def \u{{\mathbf{u}}}
\newcommand \w{{\mathbf{w}}}

\begin{document}
\title{{Definite Determinantal representations via Orthostochastic matrices}}
\author{Papri Dey}

\date{}
\maketitle
\begin{abstract}
Determinantal polynomials play a crucial role in semidefinite programming problems.  Helton-Vinnikov proved that real zero (RZ) bivariate polynomials are determinantal. However,  it leads to a challenging problem to compute such a determinantal representation. We provide a necessary and sufficient condition for the existence of definite determinantal representation of a bivariate polynomial by identifying its coefficients as  scalar products of two vectors where the scalar products are defined  by orthostochastic  matrices.  This alternative condition enables us to develop a method to compute a monic symmetric/Hermitian determinantal representations for a bivariate polynomial of degree $d$. In addition, we  propose a computational relaxation to the determinantal problem which turns into a problem of expressing the vector of coefficients of the given polynomial as convex combinations of some specified points.  We also characterize the range set of vector coefficients of a certain type of  determinantal  bivariate  polynomials.
\end{abstract}

\noin \textbf{AMS Classification (2010)}. 15A75, 15B10, 15B51, 90C22. \\
\noin \textbf{Keywords}. Semidefinite Programming, LMI Representable sets, Determinantal Polynomials, RZ Polynomials, Orthostochastic matrices, Exterior algebra. 

\section{Introduction}
One of the objectives in \textit{convex algebraic geometry} is to characterize convex  semi-algebraic sets which are definite  LMI representable sets. A set $S \subseteq \R^{n}$ is said to be \textit{LMI representable} if
\begin{equation} \label{lmiset}
S= \{ \x \in \R^{n} : A_{0} +x_{1}A_{1} +x_{2}A_{2} + \dots + x_{n}A_{n} \succeq 0 \}
\end{equation}
for some real symmetric matrices $A_{i}, i=0,\dots,n$ and $\x=(x_{1}, \dots, x_{n})^{T}$. If $A_{0}=I (\succ 0$),
the set $S$ is called a \textit{monic (definite) LMI representable} set. By $A \succ 0 (\succeq 0)$ we mean that the matrix $A$
is positive (semi)-definite. A \textit{spectrahedron} which is the feasible set of a \textit{semidefinite programming}
(SDP) problem is an another term used for a LMI representable set. %

An approximate optimal  solution of a SDP can be found by applying \textit{interior point methods} 
\cite{Nesterov}, \cite{Boydlmi}  when the SDP is strictly feasible. The assumption of strict feasibility of a 
SDP problem is equivalent to the assumption that its feasible set has nonempty interior. It is proved that if the set 
$S$ has non-empty interior, the constant coefficient matrix $A_{0}$ can be chosen to be positive definite 
[\cite{Ramana1},section $1.4$], \cite{Tim}.  So,  if the feasible set of an optimization problem is a definite 
linear matrix inequality (LMI)  representable set, the optimization problem can be transformed into a SDP 
problem \cite{Ramana1}, \cite{Helton}. The technique of converting optimization problems into semidefinite 
programming (SDP) problems arise in control theory,
signal processing  and many other areas in engineering. Now we briefly talk about the 
connection between LMI representable sets and determinantal polynomials.

A polynomial $f(\x) \in \R[\x]$ is said to be a determinantal polynomial if it can be written as
 \begin{equation} \label{detreppoly}
f(\x) = \det(A_{0}+x_{1}A_{1} +x_{2}A_{2} + \dots + x_{n}A_{n}),
\end{equation}
where coefficient matrices $A_{i}$ of linear matrix polynomial  are symmetric/Hermitian of some order greater than the degree of the polynomial and the 
constant coefficient matrix $A_{0}$ is positive definite. Then the algebraic interior associated with $f(\x)$ i.e., the 
closure of a (arcwise) connected component of $\{ \x \in \R^{n}:f(\x) > 0\}$  is a spectrahedron \cite{Helton}.
Thus one of the successful techniques  to  deal with characterizing definite LMI representable sets is to characterize
determinantal polynomials.  So, we focus on definite (monic) symmetric/Hermitian determinantal representation 
in order to accomplish the connection  between determinantal  polynomials and semidefinite programming (SDP) problems.

The determinantal  polynomials  are of special kind of polynomial, called real zero (RZ) polynomials \cite{Helton}. 
A multivariate polynomial $f(\x) \in \R[\x]$ is said to be a real zero (RZ) polynomial if the polynomial has only real 
zeros when it is restricted to any line passing through origin i.e., for any $\x \in \R^{n}$, all the roots of the 
univariate polynomial $f_{\x}(t):=f(t\cdot \x)$ are real and $f(0) \neq 0$. The polynomial $f(\x)$  is called strictly 
RZ if all these roots are distinct, for all $\x \in \R^{n}, \x \neq 0$. The homogenization of a RZ polynomial is known as hyperbolic polynomial which has a vast area of research on its own. For example, one can see \cite{Nuij}, \cite{Henrion},\cite{Brahyper}, \cite{Sturmfelsbivariate}, \cite{Vinzant}, \cite{Plaumann}, \cite{Thorstenhyp}, \cite{paprirz} from the literature.

Helton-Vinnikov have proved that a RZ bivariate polynomial $f(x_{1},x_{2})$ always admits monic Hermitian as well as 
symmetric determinantal representations of size $d$ \cite{Helton}. The homogenized version of this result is 
known as Lax conjecture \cite{Pablo3}. However, it is computationally difficult task to test whether a given polynomial is hyperbolic (with respect to a
fixed point $e$). The precise complexity is only known in some special cases (see \cite{Nikhilhyp}). 


The authors in \cite{Helton} have provided explicit expressions of the coefficient matrices of a symmetric determinantal 
representation in terms of theta functions, and the period matrix of the curve $f(x_{1},x_{2})=0$ when the curve is defined  
by a strictly RZ bivariate  polynomial $f(x_{1},x_{2})$.
Indeed, it is not easy to compute determinantal representation numerically or symbolically using this method. Later, 
the problem of  computing monic  symmetric/Hermitian determinantal representation for a strictly RZ bivariate 
polynomial has been widely studied, for example one can see  \cite{Dixon}, \cite{Sturmfelsbivariate}, \cite{Henrion}, 
\cite{Vinnikov2}.

\subsection*{The results of this paper}
In this paper, we provide a representation of the vector coefficient of mixed monomials of a determinantal bivariate polynomial as scalar product of two known vectors with different defining matrices. It is shown that defining matrices are orthostochastic matrices. This provides us a necessary and sufficient condition for the existence of a monic symmetric/Hermitian determinantal representation, see the Theorem \ref{themscalar}. 

This necessary and sufficient condition can be treated as an alternative condition for a bivariate polynomial to be a determinantal polynomial, although it's different from RZ property of a polynomial in the sense of computation. An important fact about this alternative condition is that it provides a mean to compute such a determinantal representation as opposed to RZ condition.

The proposed method can compute the eigenvalues and diagonal entries of coefficient matrices by solving roots of certain univariate polynomials which are restrictions of the given polynomial along coordinates and solving systems of linear equations respectively, see
Lemma \ref{lemmartdiag} and Corollary \ref{Corollarydiag}.
 As we need the coefficient matrices to be symmetric/Hermitian, so all the eigenvalues and diagonal entries of coefficient matrices must be real. Thus these results provide a few necessary conditions for the given polynomial of degree $d$ to be a determinantal polynomial of size $d$. 

More precisely, if such a representation exists, these two vectors in scalar product are uniquely (up to ordering) constructed from the eigenvalues of the symmetric/Hermitian coefficient matrices of a determinantal representation of the given polynomial and the defining matrices are obtained as  (complex) Hadamard product of exterior powers of an orthogonal  (unitary) matrix $V$ with themselves. In fact, these matrices are orthostochastic (resp. unistochastic) matrices corresponding to a monic symmetric (resp. Hermitian) determinantal representation. 

Moreover, we translate the determinantal representation problem into finding a suitable orthostochastic/unistochastic matrix using theory of majorization, see the Theorem \ref{themmaj}. This enables us to compute such a determinantal representation for a lower degree bivariate polynomial if it exists. The another advantage of this method is that it works even though the polynomial is not strictly RZ or the plane curve defined by the polynomial is not smooth. 

An interesting fact about this method is that it can completely characterize cubic determinantal polynomials (Section-\ref{seccubic}). Consequently, this alternative representation characterize determinantal bivariate polynomials as well as RZ bivariate polynomials.

However, we also propose a computational relaxation to determinantal problem which works well for higher degree polynomials (Section-\ref{secrelax}). The necessary and sufficient condition for the relaxation problem turns into a polytope membership problem which is computationally easier to handle. 

Explicitly, in the relaxation method one needs to check whether a point which is the vector coefficient of a bivariate polynomial can be written as convex 
combinations of some specified points. It's also shown how this method can help us to compute such a determinantal representation for bivariate polynomial by an example of quartic case in details.

In Section \ref{secrange}, we study the range set of vector coefficients of mixed monomials of certain class of determinantal bivariate polynomials and have shown that it lies inside the convex hull of some specified points. All these methods have been implemented in \textit{Macaulay2} in DeterminantalRepresentations package. 
\textit{Monic symmetric determinantal representation} is abbreviated as MSDR and  \textit{monic Hermitian determinantal representation} is abbreviated as MHDR in this paper.

\textbf{Acknowledgements.} I would like to thank my PhD supervisor Prof. Harish K. Pillai for helpful discussions on the subject of this paper. I am grateful to Prof.~J.W. Helton and Prof. Cynthia Vinzant for  useful suggestions to express the ideas of this paper properly. I would like to thank Dr. Justin Chen for helping me to implement these methods in \textit{Macaulay2}.  Much of the work on this paper  has been supported by Council of Scientific and Industrial Research (CSIR), India while the author was doing her doctoral work in IIT Bombay. The author also gratefully acknowledges support through the
\textit{Max Planck Institute for Mathematics in the Sciences} in Leipzig, Germany and the Institute
for Computational and Experimental Research in Mathematics, Brown University, USA.
\section{Determinantal Polynomials}
In this section, we convert the determinantal representations of bivariate polynomials into another problem which is concerned to find a suitable orthostochastic or unistochastic matrix corresponding to MSDR or MHDR respectively.
\subsection{Determining the Eigenvalues of Coefficient Matrices:}
First we notice some facts about determinantal multivariate polynomials. 
Since the coefficient matrices of a determinantal polynomial are symmetric (Hermitian), therefore by the
spectral theorem of a symmetric (Hermitian) matrix there exist 
a suitable orthogonal (unitary) matrix $U$ such that one of the coefficient matrices becomes diagonal. Without loss of generality, it is enough to consider coefficient matrix associated to $x_{1}$ as a diagonal matrix $D_{1}$ and obtain an MSDR (MHDR) of the following form
\begin{equation}
f(\x)=
\det(I+x_{1}D_{1}+x_{2}A_{2}\dots+x_{n}A_{n})
\end{equation}

We explain a technique to determine the eigenvalues of the coefficient matrices $D_{1}$ and $A_{i}, 2\leq i \leq n$. 
We take restrictions of the given multivariate polynomial $f(\x)$ along each  $x_{i}, i=1, \dots,n$ that means we  restrict the polynomial along one variable at a time by making the rest of the variables zero and generate $n$ univariate polynomials $f_{x_{i}}=f(0,\dots,x_{i},\dots,0)$.

It is known that if a multivariate polynomial $f(\x)$ admits an MSDR (MHDR), it is a RZ polynomial. By recalling the definition of RZ polynomial, we know that for any $\x \in \R^{n}$, RZ polynomial $f(\x)$ when restricted along any line passing through origin, i.e., the univariate polynomial $f_{\x}(t)$ has only real zeros. So when a RZ polynomial $f(\x)$ restricted along $x_{i}, i=1,\dots,n$, each of them has only real zeros, i.e., all univariate polynomials $f_{x_{i}}$ in $x_{i}$ have only real zeros.

As a consequence of this result we have a necessary condition for the existence of an MSDR (MHDR) of size equal to the degree of the polynomial for a multivariate polynomial of any degree.
\begin{lemma}\label{existence of diagonal 2}
If a multivariate polynomial $f(\x) \in \R[\x]$ of degree $d$ has an MSDR (MHDR) of size $d$, then all the roots of $f_{x_{i}}$ are real for all $i=1,\dots,n$.
\end{lemma}

Moreover, the eigenvalues of matrices $A_{i}$ can be found from the roots of $f_{x_{i}}$ for all $i=1, \dots, n$ by the following Lemma which is proved in \cite{Paprimulti}. For the sake of completeness we include the proof here. 
\begin{lemma} \label{lemmartdiag}
The eigenvalues of coefficient matrices $A_{i}$ are the negative reciprocal of the roots of univariate polynomials $f_{x_{i}}$ for all $i=1, \dots,n$.
\end{lemma}
\pf The eigenvalues of coefficient matrices $A_{i}$ are all real as the coefficient matrices are either symmetric or Hermitian. By Lemma \ref{existence of diagonal 2}, all the roots of $f_{x_{i}}$ are real for all $i=1,\dots,n$. If a univariate polynomial $f(x)=\det(xI+A)$ of degree $d$ has only real zeros, so is the reversed polynomial $\tilde{f}(x):=x^{d}f(1/x)$. 

On the other hand, $\det(tI+A_{i})=t^{d}\tilde{f}_{e_{i}}(e_{i}/t)$ at $\x=e_{i}$, where $e_{i}$ denotes the standard basis vector in $\R^{n}$. 
Therefore, there is a one to one correspondence between the roots of $f_{x_{i}}$ and the non-zero eigenvalues of $A_{i}$ and here the map is $t \mapsto -1/t$.  \qed
\begin{remark} \rm{
 Eigenvalues of $A_{i},i=1,2$ are unique up to ordering. Without loss of generality arrange them in descending order.
 }
\end{remark}
Now we state the following theorem which talks about the relations between the coefficients  of a determinantal polynomial and the coefficient matrices of its determinantal representation by using the notion of  generalized mixed discriminant of coefficient matrices \cite{Paprimulti}. 
\begin{theorem} (Generalized Mixed Discriminant Theorem) \label{themgmd}
The coefficients of a determinantal polynomial $f(\x)$ of degree $d$ are uniquely determined by the generalized mixed discriminants of the coefficient matrices $A_{i}$ as follows. If the degree of a monomial $x_{1}^{k_{1}}x_{2}^{k_{2}} \dots x_{n}^{k_{n}}$ is $k$ (i.e., $k_{1}+k_{2}+\dots+k_{n}=k) \leq d$, then the coefficient 
of ($x_{1}^{k_{1}} \dots x_{n}^{k_{n}})$ is given by 
\begin{equation*}
\widehat{D}(\underbrace{A_{1}, \dots , A_{1}}_{k_{1}}, \underbrace{A_{2}, \dots, A_{2}}_{k_{2}},\dots,\underbrace{A_{n},\dots,A_{n}}_{k_{n}}).
\end{equation*}
\end{theorem}
Then using the generalized mixed discriminant Theorem we compute the diagonal entries of $A_{i},i=1,2$.
Let
\begin{eqnarray*} \label{eqscalar}
&\eig(A_{1}):=\u_{1}:=\bmatrix{d_{1}&d_{2}& \dots & d_{d}}^{T}, \eig(A_{2}):=\w_{1}:=\bmatrix{r_{1}&r_{2}& \dots & r_{d}}^{T}, \\ &\y_{i}:=\diag(A_{i}),i=1,2
\end{eqnarray*}
By the Theorem \ref{themgmd} the diagonal entries of coefficient matrices $A_{i},i=1,2$ can be determined by solving systems of linear equations of the form $G_{i} \y_{i} =\z_{i}$, where   
\begin{equation} \label{eqdiag}
\z_{i}=\bmatrix{\mbox{coeff of} \ x_{i} \\ \mbox{coeff of} \ x_{i}x_{j}\\ \vdots  \\ \mbox{coeff of}\  x_{i}x_{j}^{d-1}}_{j \neq i}, G_{1}=\bmatrix{1 & 1 & \dots & 1\\ \sum_{i=2 }^{d}r_{i}& \sum_{i=1, i \neq 2}^{d}r_{i}& \dots &  \sum_{i=1}^{d-1}r_{i}\\ \sum_{i_{k},i_{l}\neq 1,i_{k} <i_{ l}}r_{i_{k}}r_{i_{l}} & \dots & \dots & \sum_{i_{k},i_{l}\neq d,i_{k} <i_{ l}}r_{i_{k}}r_{i_{l}}  \\ \vdots & \vdots & \vdots & \vdots \\ r_{2}r_{3}\dots r_{d} & r_{1}r_{3}\dots r_{
d} & \dots & r_{1}\dots r_{d-1}}.
 \end{equation}
Note that $G_{1}$ involves only the eigenvalues of $A_{2}$. Similarly, matrix $G_{2}$ which is defined by replacing $r_{i}$ with $d_{i}$ in equation (\ref{eqdiag}) involves only the eigenvalues of 
$A_{1}$.

As $\z_{1}$ is the vector of coefficients of monomials of the form $x_{1}x_{2}^{\alpha_{2}}, 0 \leq \alpha_{2} \leq d-1$, so 
the relations are linear in terms of the entries of diagonal entries of $A_{1}$ by the Theorem \ref{themgmd}. In fact, this makes it possible to compute diagonal entries of $A_{1}$ by solving a system of linear equations. Similarly, one can compute the diagonal entries of $A_{2}$ by considering the vector $\z_{2}$ associated with monomials $x_{1}^{\alpha_{1}}x_{2}, 0 \leq \alpha_{1} \leq d-1$ .

Note that $\det(G_{1})=\prod_{i < j}(r_{i}-r_{j})$ and $\det(G_{2})=\prod_{i < j}(d_{i}-d_{j})$. Therefore, the  diagonal entries of $A_{1},A_{2}$ are uniquely determined up to ordering  if all the eigenvalues of coefficient  matrices $A_{2},  A_{1}$ are distinct respectively.
\begin{corollary} \label{Corollarydiag}
The diagonal entries of coefficient matrices $(A_{1},A_{2})$ of a determinantal polynomial $f(\x)=\det(I+x_{1}A_{1}+x_{2}A_{2})$ can be determined (uniquely up to ordering) by solving the systems of linear equations defined in equation (\ref{eqdiag}) (provided $G_{i}$ is invertible).  
\end{corollary}

\subsection{Connections With Orthostochastic (Unistochastic) Matrices} 
In this section, we provide a necessary and sufficient condition for the existence of an MSDR (MHDR) of size $d$ of a bivariate polynomial of degree $d$ by representing its coefficients as scalar products of two vectors defined by different matrices. Note that this technique works well even if the polynomial is not a strictly RZ  polynomial that means repeated eigenvalues of coefficient matrices are allowed.

At first, we briefly recall some basic definitions and facts that will be used in sequel.   A \textbf{doubly stochastic} matrix is a square matrix whose entries are nonnegative and the sum of the elements in each row and each column is unity. An \textbf{othostochastic} matrix is a doubly stochastic matrix whose entries are the squares of the entries of some orthogonal matrix. A \textbf{unistochastic} matrix is a doubly stochastic matrix whose entries are the squares of the absolute values of the entries of some unitary matrix.

A bilinear form on $\R^{n}$ is a map from $\R^{n} \times \R^{n}$ to $\R$ defined by
$(\x, \y) \mapsto \langle \x, \y \rangle_{Q}=\x^{T}Q\y$ where the  matrix $Q \in \R^{n \times n}$ is
associated with the bilinear form for all $\x, \y \in \R^{n}$. The matrix $Q$ is known as the defining
matrix of the bilinear form. The form is said to be non degenerate when $Q$ is nonsingular.
This is also known as scalar product.

\textbf{Constructions of Vectors in Scalar Products}:
Let $\N_{k}=\{\delta:=(i_{1},\dots,i_{j},\dots,i_{k}) \in \N^{k}: i_{1} < \dots < i_{k}, 1\leq k \leq d\}$ denotes the $k$ ordered index set.

The $\delta$-th components of $\u_{k}$ and $\w_{k}$ are the the product of $k$ ordered (with $i_{1} < \dots < i_{k})$  eigenvalues of $A_{1}$ and $A_{2}$ respectively, i.e.,
\begin{equation}  \label{eqscalar1}
 \u_{k}(\delta)=\prod_{i \in \delta}d_{i}, \ \  \w_{k}(\delta)=\prod_{i \in \delta}r_{i},  \delta \in \N_{k}.
\end{equation}
Observe that the multi-indexed vectors $\u_{k}$ and $\w_{k}$ can be uniquely determined up to ordering of eigenvalues of $A_{i},i=1,2$. Now we define another set of multi-indexed vectors.  The $\delta'$-th component of $\u_{k,k'}^{c}$ is the sum of the product of all possible combinations of $k$ ordered eigenvalues of $A_{1}$ except those combinations which involve at least one of the eigenvalues of $\delta'$-th component of $\u_{k'}$, and the $\delta'$-th component of $\w_{k,k'}^{c}$ is the sum of the product of all possible combinations of $k$ ordered eigenvalues of $A_{2}$ except those combinations which involve at least one of the eigenvalues of $\delta'$-th component of $\w_{k'}$, i.e.,
\begin{equation} \label{eqscalar2}
\u_{k,k'}^{c}(\delta')=\sum_{(\delta,\delta') \in \N_{k} \times \N_{k'}, \delta \cap \delta' =\emptyset} \prod_{i \in \delta} d_{i},\ \ \w_{k,k'}^{c}(\delta')=\sum_{(\delta,\delta') \in \N_{k} \times \N_{k'}, \delta \cap \delta' =\emptyset} \prod_{i \in \delta} r_{i}
\end{equation}
Note that the degree of each component of $\u_{k,k'}^{c}$ depends on the degree of each component of $\u_{k}$, the support and the order  (i.e., number of components) of $\u_{k,k'}^{c}$ depend on the vector $\u_{k'}$.

\textbf{Construction of Defining Matrices in Scalar Products}:
We use the notion of $k$ -th exterior power of an orthogonal matrix. In order to explain what is meant by $k$ -th exterior power of an orthogonal matrix which will be used in sequel we need to discuss some preliminaries \cite{Sergeiexterior}, \cite{Taylor}


If $\{e_{1},e_{2},\dots,e_{d}\}$ is the standard basis of the vector space $\E^{d}$ over the field $\K(\R$ or $\C$), then the set of $k$-vectors $e_{I}:=e_{i_{1}} \wedge \dots \wedge e_{i_{k}}$ where $I=\{(i_{1},i_{2},\dots,i_{k}): 1 \leq i_{1} < \dots < i_{k} \leq d\}$ is a basis of the $k$-th exterior
power of $\E^{d}$, denoted by $\wedge^{k}\E^{d}$ for $k=1, \dots,d$ and thus any element in $\wedge^{k}\E^{d}$ which is nothing but a $k$ vector (sum of $k$ blades) can be uniquely written as
\begin{equation*}
\a=\sum_{I \subseteq \{1,\dots,d\}, |I|=k}a_{I}e_{I}=\sum_{1 \leq i_{1} < \dots < i_{k} \leq d}a_{i_{1}\dots i_{k}}(e_{i_{1}}\wedge \dots \wedge e_{i_{k}}).
\end{equation*}
For $k=1,\dots,d$, the norm of $\a$ is $|\a|= \sqrt{\sum_{i_{1}< \dots < i_{k}} (a_{i_{1}\dots i_{k}})^{2}}$ and we set $|\a|=0, k >d$.
The homogeneous coordinates $a_{I}$s are known as Pl\"{u}cker coordinates on $\P(\wedge^{k}\E^{d})$ 
 associated with the ordered basis $(e_{1},\dots,e_{d})$ of $\E^{d}$. Naturally, if $a_{1},\dots, a_{d} \in \E^{d}$, then
\begin{equation*}
a_{1} \wedge \dots \wedge a_{d}=\det(a_{1}^{T},\dots,a_{d}^{T})e_{1} \wedge \dots \wedge e_{d}
\end{equation*}

The collection of the spaces $\wedge^{k}(\E^{d})$, for $k = 0, 1, 2,\dots$ , together with the
exterior product or operation $\wedge$ is called the \textit{Exterior algebra} or \textit{Grassmannian algebra}
on $\E^{d}$. So, we have $\wedge(\E^{d})= \bigoplus_{k=0}^{\infty} \wedge^{k}(\E^{d})=\bigoplus_{k=0}^{d} \wedge^{k}(\E^{d})$, as $\wedge^{k}(\E^{d})=0,$ if $k > d$.

The set of all $d \times d$ orthogonal  (resp.unitary)  matrices is the orthogonal  (resp.unitary)  group $O(d) (U(d))$. The group $O(d)$  (resp. $U(d))$ can be
identified by the set of its $d$-tuple of ordered column vectors. Mathematically,
\begin{align*}
O(d)&=\{V:=(v_{1},\dots,v_{d}) \in \R^{d} \times \dots \times \R^{d}, v_{i}^{T}v_{j}=\delta_{ij} \} \\
U(d)&=\{U:=(u_{1},\dots,u_{d}) \in \C^{d} \times \dots \times \C^{d}, u_{i}^{\ast}u_{j}=\delta_{ij} \}
\end{align*}
where $\delta_{ij}$ is the Kronecker delta. In our context, the $k$-th exterior power of a matrix $V:=(v_{1},\dots,v_{d})$ can be identified by the set of its $d \choose k$-tuple of ordered column vectors, denoted by $V^{\wedge^{k}}$ and each column of $V^{\wedge^{k}}$ is equal to $(v_{i_{1}} \wedge \dots \wedge v_{i_{k}})$ 
where $\{(i_{1},\dots,i_{k}):i_{1} < \dots < i_{k}, i_{j} \in \{1,\dots,d\}\}$ is an ordered $k$ tuple set.

The (complex) Hadamard product of $k$-th exterior power of  an orthogonal (unitary) matrix
$V$ with itself (its complex conjugate) is denoted by
$Q^{\wedge^{k}}:=(V^{\wedge^{k}}\odot V^{\wedge^{k}})$ for all $k=2,\dots,n$.
Thus the $ij$ element of the matrix $Q^{\wedge^{k}}$ is defined as the square of $k \times k$ minor of
the matrix $V^{\wedge^{k}}$ and $k \times k$ minors are the determinants of matrices of order $k$
constructed by choosing rows corresponding to $i$th component of $\u_{k}$ and columns corresponding
to $j$th component of $\w_{k}$. In particular, $Q=(v_{ij}^{2})$, if $V=(v_{ij})$.

Note that the set of columns of an orthogonal matrix $V \in O(d)$ which forms an orthonormal basis of the Euclidean space $\R^{d}$ generates an orthonormal basis  by identifying the columns of  $V^{\wedge^{k}}$ of the vector space $\wedge^{k}\R^{d}$  \cite{Vincentext}. Similarly, the set of columns of a unitary matrix $U \in U(d)$ which forms an orthonormal basis of the vector space $\C^{d}$ generates an orthonormal basis  by identifying the columns of $U^{\wedge^{k}}$ of the vector space $\wedge^{k}\C^{d}$ \cite{Vincentext}.

Using these facts we  prove the following results.
\begin{lemma}\label{orthostochasticmatrix}
The Hadamard product of $k$-th exterior power matrix of the orthogonal matrix $V$ of order $d$ with itself,
denoted as $Q^{\wedge^{k}}:=V^{\wedge^{k}} \odot V^{\wedge^{k}}$  is an orthostochastic matrix for all $1 \leq k \leq d$.
\end{lemma}
\pf From the construction of the matrix $Q^{\wedge^{k}}$, it is clear that each entry is
a square of some $k \times k$ minors where $1 \leq k \leq d$. So, $Q^{\wedge^{k}}(ij) \geq 0$. Let $V=\bmatrix{\vdots & v_{j}  & \vdots}=\bmatrix{\dots \\ w_{j} \\ \dots} \in O(d);j=i \dots,d$ be an orthogonal matrix. As the matrix $V$ is an orthogonal matrix, therefore $\{v_{1},\dots,v_{d}\}$ and $\{w_{1}, \dots, w_{d}\}$ are sets of orthonormal vectors.
Therefore, the set of $k$ vectors of the form
$v_{i_{1}}\wedge \dots \wedge v_{i_{k}}$ where $i_{1} < \dots < i_{k}$, is a basis for
$\wedge^{k}\R^{d}$ and satisfies
\begin{equation*}
\langle v_{i_{1}} \wedge \dots \wedge v_{i_{k}}, v_{j_{1}} \wedge \dots \wedge v_{j_{k}} \rangle=
\left\{
\begin{array}{lcl}
1 \ &\mbox{if} \ (i_{1},\dots,i_{k})=(j_{1},\dots,j_{k}) \\
0 \ &\mbox{otherwise}.
\end{array}
\right.
\end{equation*} \\
Note that the set $\{v_{i_{1}} \wedge \dots \wedge v_{i_{k}}, i_{1} < \dots < i_{k} \}$ is a
collection of $d \choose k$ orthonormal $k$-vectors and each of these orthonormal $k$-vectors represents a column of matrix $V^{\wedge^{k}}$ . Thus 
\begin{equation*}
 Q^{\wedge^{k}}=\bmatrix{\vdots & v_{i_{1}} \wedge \dots \wedge v_{i_{k}} & \vdots} \odot \bmatrix{\vdots & v_{i_{1}} \wedge \dots \wedge v_{i_{k}} & \vdots} =\bmatrix{\dots \\ w_{i_{1}} \wedge \dots \wedge w _{i_{k}}\\ \dots} \odot \bmatrix{\dots \\ w_{i_{1}} \wedge \dots \wedge w _{i_{k}}\\ \dots}.
 \end{equation*} So, each row sums and column sums of these matrices are actually $|v_{i_{1}}\wedge \dots \wedge v_{i_{k}}|^{2}=1$
and $|w_{i_{1}}\wedge \dots \wedge w_{i_{k}}|^{2}=1$
respectively. Therefore,
$\sum_{i}Q^{\wedge^{k}}(ij)=\sum_{j}Q^{\wedge^{k}}(ij)=1$.
So, the matrix $Q^{\wedge^{k}}$ is an orthostochastic matrix for all $1 \leq k \leq d$. \qed 

Similarly, we conclude that
\begin{corollary}\label{corunistochastic}
The complex Hadamard product of $k$-th exterior power matrix of the unitary matrix $U$ of order $d$ with itself,
denoted by $Q^{\wedge^{k}}:=U^{\wedge^{k}} \odot (U^{\wedge^{k}})^{\ast}$ is unistochastic matrix for all $1 \leq k \leq d$.
\end{corollary}
\subsection{Necessary and Sufficient Condition:}
In this subsection, we propose a representation of coefficients of mixed monomials of a bivariate determinantal polynomial
which provides a necessary and sufficient condition for the existence of MSDRs (MHDRs) of size $d$  for a bivariate polynomial of degree $d$.

We have shown that the eigenvalues of coefficient matrices are uniquely  determined by using the coefficients of
monomials $x_{i}^{\alpha_{1}}, \alpha_{1} \in \{0,\dots,d\},i=1,2$. Thus these coefficients of
monomials in $x_{i}, i=1,2$ of a bivariate polynomial can be expressed in terms of the eigenvalues of the
corresponding coefficient matrices and they are independent of the choice of orthogonal (unitary) matrix.  

The choice of orthogonal (unitary) matrix affects only on the vector coefficient of mixed monomials of a determinantal bivariate polynomial. By mixed monomials we mean to specify the monomials which are consisting of all the variables with at least of degree one or in other words,   each variable should appear in those monomials. So, in order to find a suitable orthogonal (unitary) matrix it is enough to study the behaviour of the vector coefficient of mixed monomials of a bivariate polynomial.

Consider the bivariate polynomial $f(\x)=f_{d0}x_{1}^{d}+\dots+f_{10}x_{1}+f_{0d}x_{2}^{d}+\dots+f_{01}x_{2}+\widetilde{f}(\ x)+1$ of degree $d$, where $\widetilde{f}(\x)=\sum_{\alpha_{1},\alpha_{2}\geq 1}^{d-1} f_{\alpha_{1}\alpha_{2}} x_{1}^{\alpha_{1}}x_{2}^{\alpha_{2}}$,
$|f_{\alpha_{1}\alpha_{2}}|=$ $n+d-1 \choose d$ $-n$.
\begin{theorem} \label{themscalar}
A bivariate polynomial $f(\x)$ of degree $d$ has an MSDR  of size $d$ if and only if there exists an orthostochastic matrix $Q:=V \odot V$ such that 
\begin{enumerate}
\item Case I: $\alpha_{1} \geq \alpha_{2}$.
 The coefficient $f_{\alpha_{1}\alpha_{2}}=$
 ${\u_{\alpha_{1}, \alpha_{2}}^{c}}^{T}Q^{\wedge^{\alpha_{2}}} \w_{\alpha_{2}}$
\item Case II: $\alpha_{1} \leq \alpha_{2}$. The coefficient $f_{\alpha_{1}\alpha_{2}}=$ 
${\w_{\alpha_{2}, \alpha_{1}}^{c}}^{T}(Q^{\wedge^{\alpha_{1}}})^{T} \u_{\alpha_{1}}$ 
\end{enumerate}
where vectors $\u_{\alpha_{1}}, \u_{\alpha_{1},\alpha_{2}}^{c}, \w_{\alpha_{1}},\w_{\alpha_{2},\alpha_{1}}^{c}$ are determined by the eigenvalues of coefficient matrices as defined in equations (\ref{eqscalar1}), (\ref{eqscalar2}) and $Q^{\wedge^{k}}$ denotes the Hadamard product of $k$ th exterior power of an orthogonal matrix $V$ with itself. 
\end{theorem}
\pf A bivariate polynomial $f(\x)$ of degree $d$ has an MSDR   of size $d$ if and only if there exists an orthogonal matrix $V$ such that
\begin{equation} \label{eqdet}
f(\x)=\det(I+x_{1}D_{1}+x_{2}VD_{2}V^{T})
\end{equation}
By Lemma \ref{lemmartdiag} the entries of diagonal matrices $D_{1}$ and $D_{2}$ are uniquely determined by the vector coefficients of monomials associated with univariate polynomials $f_{x_{1}}$ and $f_{x_{2}}$ respectively. This implies that the eigenvalues of coefficient matrices are uniquely determined by these vector coefficients. Thus, the vectors $\u_{1}, \w_{1}$ defined in equation (\ref{eqscalar}) are uniquely determined up to descending (ascending) order. The ordering of vectors $\u_{1},\w_{1}$ determines the vectors $\u_{\alpha_{1}}, \u_{\alpha_{1},\alpha_{2}}^{c}, \w_{\alpha_{1}},\w_{\alpha_{2},\alpha_{1}}^{c}$ uniquely up to (graded lexicographic) monomial order by using equation (\ref{eqscalar1}) and equation (\ref{eqscalar2}). Observe that the choice of orthogonal matrix $V$ affects only the vector of coefficients of mixed monomials of the given polynomial $f(x_{1}, x_{2})$, i.e., the part $\tilde{f}(\x)$ of given $f(\x)$. So, it is enough to look into representation of the vector coefficient of mixed monomials of the given polynomial $f(x_{1}, x_{2})$. The Theorem \ref{themgmd} reveals that the analytic expressions of the vector coefficient of mixed monomials $x_{1}^{\alpha_{1}}x_{2}^{\alpha_{2}} (\alpha_{1} \geq \alpha_{2})$ of a bivariate polynomial involve only the diagonal entries of the matrices $VD_{2}V^{T}$ and  $V^{\wedge^{k}}D_{2}^{\wedge^{k}}(V^{\wedge^{k}})^{T}$ where $V^{\wedge^{k}}$ denotes $k$-th exterior power of a matrix $V$. This result can be visualized from the equation (\ref{eqdet}) by simple calculations. Similarly, due to symmetry between the coefficients, the analytic expressions of the vector coefficient of mixed monomials $x_{1}^{\alpha_{1}}x_{2}^{\alpha_{2}} (\alpha_{1} \leq \alpha_{2})$ of a bivariate polynomial involve only the diagonal entries of the matrices $V^{T}D_{1}V$ and  $(V^{\wedge^{k}})^{T}D_{1}^{\wedge^{k}}V^{\wedge^{k}}$. On the other hand, since $(VD_{2}V^{T})_{ii}=\sum_{j=1}^{d}v_{ij}r_{j}v_{ij}=[(V \odot V)\w_{1}]_{i}$, where $r_{j}$ are the diagonal entries of $D_{2}$, so the $i$-th diagonal entry of the matrix $V^{\wedge^{k}}D_{2}^{\wedge^{k}}(V^{\wedge^{k}})^{T}$ coincides with the $i$-th entry of the vector $[(V^{\wedge^{k}} \odot V^{\wedge^{k}})\w_{k}]_{i}, i=1, \dots,m$. Note that $D_{2}^{\wedge^{k}}$ is a diagonal matrix whose $i$-th diagonal entry is $i$-th component of vector $\w_{k}$, and $\odot$ denotes the Hadamard product. Hence the proof. \qed

Consequently, applying the same logic the following result provides a necessary and sufficient condition for the existence of an MHDR of a bivariate polynomial of degree $d$. Note that orthogonal matrix $V$ would be replaced by unitary matrix $U$ and the $i$-th diagonal entry of the unitary matrix $U^{\wedge^{k}}D_{2}^{\wedge^{k}}(U^{\wedge^{k}})^{\ast}$ coincides with the $i$-th entry of the vector $[(U^{\wedge^{k}} \odot (U^{\wedge^{k}})^{\ast})\w_{k}]_{i}, i=1, \dots,m$ where $\odot$ denotes the complex Hadamard product. 

\begin{remark} \rm{
A bivariate polynomial $f(\x)$ of degree $d$ has an MHDR  of size $d$ if and only if there exists a unistochastic matrix $Q:=U \odot U^{\ast}$ 
\begin{enumerate}
\item Case I: $\alpha_{1} \geq \alpha_{2}$.
 The coefficient $f_{\alpha_{1}\alpha_{2}}=$ 
 ${\u_{\alpha_{1},\alpha_{2}}^{c}}^{T}Q^{\wedge^{\alpha_{2}}} \w_{\alpha_{2}}$.
\item Case II: $\alpha_{1} \leq \alpha_{2}$. The coefficient $f_{\alpha_{1}\alpha_{2}}=$ 
${\w_{\alpha_{2},\alpha_{1}}^{c}}^{T}(Q^{\wedge^{\alpha_{1}}})^{T} \u_{\alpha_{1}}$ 
\end{enumerate}
where vectors $\u_{\alpha_{1}}, \u_{\alpha_{1},\alpha_{2}}^{c}, \w_{\alpha_{1}},\w_{\alpha_{2},\alpha_{1}}^{c}$ are determined by the eigenvalues of coefficient matrices as defined in equation (\ref{eqscalar1}) and equation (\ref{eqscalar2}) and $Q^{\wedge^{k}}$ denotes the complex Hadamard product of $k$ th exterior power of a unitary matrix $U$ with itself. 
}
\end{remark}

Therefore, using the Lemma \ref{orthostochasticmatrix} we conclude that
\begin{corollary}
 The vector coefficient of mixed monomials of a determinantal polynomial $f(\x)$ can be expressed as scalar product of two vectors with orthostochastic or unistochastic defining matrices.
\end{corollary}

So, in order to determine MSDR (MHDR) our aim is to find an orthostochastic matrix (a unistochastic matrix) $Q$ which satisfies all the scalar product expression for a given bivariate polynomial. Interestingly, this issue is highly related to a well established field known as theory of majorization and also connected to the inverse eigenvalue problem.
In fact,  the following theorem \cite{Horn} provides a necessary and sufficient condition for the existence of such an orthostochastic (a unistochastic) matrix for a pair of majorized vectors in the field of majorization theory.

\begin{theorem}  \label{orthoconv2} \cite{Horn}
Let $\x,\y \in \R^{d}$. Then the following statements are equivalent:
\begin{enumerate}
\item $\y$ is majorized by $\x$, denoted by $\y \prec \x$. By definition of majorization the following conditions are satisfied. $\max _{\sigma \in S^{d}_{k}} \sum _{i=1}^{k} y_{\sigma_{i}}  \leq \max _{\sigma \in S^{d}_{k}} \sum _{i=1}^{k} x_{\sigma_{i}}, 1 \leq k \leq d,$ and   $\sum_{i=1}^{d} y_{i} =\sum _{i=1}^{d} x_{i}$, where $S^{d}_{k}$ is the set of all $k$ termed sequences $\sigma$ of integers such that $1 \leq
    \sigma_{1} < \dots < \sigma_{k} \leq d$. 
\item $\y \in C(\x)$, where $C(\x)$ is the convex hull of all the points $(\x_{\alpha_{1}},\dots,\x_{\alpha_{d}}), \alpha$ varying over all permutations
of $(1,\dots,d)$.

\item $\y=Q\x$ for some orthostochastic matrix $Q$.
\end{enumerate}
\end{theorem}
Horn \cite{Horn} proved that a Hermitian matrix $H$ with eigenvalues $\x$ and diagonal entries $\y$ exists if and only if $\x$ majorizes $\y$. Later due to Horn and Mirsky, it is proved that there exists a symmetric matrix with eigenvalues $\x$ and diagonal entries $\y$ if and only if $\x$ majorizes $\y$ \cite{Marshalinequality}.

Let $\x,\y \in \R^{d}$ be arranged in descending order i.e., $x_{1} \geq \dots \geq x_{d}$ and $y_{1} \geq \dots \geq y_{d}$. So, $\y \prec \x$ on $\R^{d}$ if and only if there
exists a Hermitian as well as a symmetric matrix with diagonal elements $y_{1},\dots, y_{d}$
and eigenvalues $x_{1},\dots, x_{d}$.

Based on these results we have a necessary condition which involves the eigenvalues and diagonal entries of coefficient matrix associated with variable $x_{2}$.

Due to symmetry between the coefficients of bivariate polynomial we could choose the coefficient matrix associated with $x_{2}$ as diagonal matrix and then by using the same terminology we could derive another necessary condition which
involves the eigenvalues and diagonal entries of coefficient matrix associated with variable $x_{1}$. 

Therefore, we provide two more necessary conditions at this stage.
\begin{proposition} \label{propnec}
 If a polynomial $f(\x) \in \R[\x]$ is determinantal, $\diag(A_{i}) \prec \eig(A_{i})$ for all $i=1,2$.
\end{proposition}


\begin{remark}\label{remarkindneccond}
\rm{ These two necessary conditions are independent and they are not sufficient condition.  It is shown in the example \ref{examcubicdet}.}
\end{remark}
Consider the set of all $d \times d$ doubly stochastic matrices, known as Birkhoff Polytope $\Omega_{d}$.
Let $\Omega_{d}(\y \prec \x)=\{Q \in \Omega_{n}; \y= Q \x, \x, \y \in \R^{d}\}$.
Then the set $\Omega_{d}(y \prec x)$ is a nonempty, convex polytope and a subpolytope of $\Omega_{d}$
[\cite{Brualdi}, Chapter-$9$]. This set is known as \textit{doubly stochastic polytope of the majorization}
$\y \prec \x$. Thus, we provide a necessary and sufficient for the existence of an MSDR (MHDR) of size $d$ for a bivariate polynomial.
\begin{theorem}\label{themmaj}
A bivariate polynomial $f(x_{1},x_{2})$ of degree $d$ admits an MSDR (MHDR) of size $d$ if and only if there exists an orthostochastic (a unistochastic) matrix $Q$ such that
$\diag(A_{1})=Q^{T} \eig(A_{1})$ and $\diag(A_{2})=Q ~ \eig(A_{2})$ and 
for all $\alpha_{1},\alpha_{2} \in \{2,\dots,d-2\}$.
\begin{enumerate}
\item Case I: $\alpha_{1} \geq \alpha_{2}$.
 The coefficient $f_{\alpha_{1}\alpha_{2}}=$
 ${\u_{\alpha_{1},\alpha_{2}}^{c}}^{T}Q^{\wedge^{\alpha_{2}}} \w_{\alpha_{2}}$
\item Case II: $\alpha_{1} \leq \alpha_{2}$. The coefficient $f_{\alpha_{1}\alpha_{2}}=$ 
${\w_{\alpha_{2}, \alpha_{1}}^{c}}^{T}(Q^{\wedge^{\alpha_{1}}})^{T} \u_{\alpha_{1}}$ 
\end{enumerate}
where vectors $\u_{\alpha_{1}}, \u_{\alpha_{1},\alpha_{2}}^{c}, \w_{\alpha_{1}},\w_{\alpha_{2},\alpha_{1}}^{c}$ are determined by the eigenvalues of coefficient matrices as defined in equation (\ref{eqscalar1}) and equation (\ref{eqscalar2}) and $Q^{\wedge^{k}}$ denotes the Hadamard product of $k$ th exterior power of an orthogonal matrix $V$ with itself.
\end{theorem}
\pf By the Theorem \ref{themscalar} a bivariate polynomial of degree $d$ admits an MSDR (MHDR) of size $d$ if and only if
there exists an orthostochastic (unistochastic) matrix such that the vector coefficient of mixed monomials are satisfied. So, the necessary part follows from the Theorem \ref{themscalar} and Proposition \ref{propnec}.

The idea of the sufficient part of this theorem is to minimize the number of mixed monomials conditions. If there exists an orthostochastic (a unistochastic) matrix $Q$ such that $Q^{T} \in \Omega_{d}(\diag(A_{1})) \prec \eig(A_{1})$ and $Q \in \Omega_{d}(\diag(A_{2})) \prec \eig(A_{2})$, the vector coefficient of mixed monomials of the form $x_{1}^{\alpha_{1}}x_{2}^{\alpha_{2}}$ in which at least one of $\alpha_{1},\alpha_{2}$ being
equal to one are already satisfied. So, the remaining monomials are the monomials $x_{1}^{\alpha_{1}}x_{2}^{\alpha_{2}}$ where $ \alpha_{1},\alpha_{2} \in \{2,\dots,d-2\}$. Hence we conclude the claim by the Theorem \ref{themscalar}.  \qed
\subsection*{Cubic bivariate Polynomials:} \label{seccubic}
We study the cubic bivariate determinantal polynomials in details. Consider the cubic bivariate polynomial 
 \begin{equation} \label{cubicbiveq1}
 f(x_{1},x_{2})=
 f_{30}x_{1}^{3}+f_{03}x_{2}^{3}+f_{21}x_{1}^{2}x_{2}+f_{12}x_{1}x_{2}^{2}+f_{20}x_{1}^{2}+f_{02}x_{2}^{2}+f_{11}x_{1}x_{2}+f_{10}x_{1}+f_{01}x_{2}+1.
 \end{equation}
Suppose the polynomial $f(\x)$ has an MSDR (MHDR) i.e., $f(\x)=\det(I+x_{1}A_{1}+x_{2}A_{2})$ where $A_{1},A_{2}$ are symmetrices of order $3$.

Note that the case of cubic bivariate determinantal polynomial is easier as the remaining coefficients due to the mixed monomials
$x_{1}^{\alpha_{1}}x_{2}^{\alpha_{2}}, \alpha_{1},\alpha_{2} \in \{2,\dots,d-2\}$ don't appear here. In fact, we can completely characterize cubic bivaraite determinantal polynomial of size $3$ using the method of finding a suitable orthostochastic or unistochastic matrix. Moreover, the number of orthostochastic matrices corresponds to the number of orthogonally non-equivalent orbits of determinantal representation.  

We can find $\eig(A_{i}),\diag(A_{i}), i=1,2$ and check whether the two pairs of majorization criteria $\diag(A_{i}) \prec \eig(A_{i})$ hold. Note that these are necessary conditions for a polynomial to be determinantal.  Although by the Theorem \ref{themscalar} we know that cubic polynomial $f(x_{1},x_{2})$ is determinantal if and only if the following conditions are satisfied. 

There exists an orthostochastic (a unistochastic) matrix $Q$ such that
\begin{equation}\label{coeffeq2}
f_{11}= {\u_{1,1}^{c}}^{T}Q \w_{1}={\w_{1,1}^{c}}^{T}Q^{T}\u_{1}, f_{21}= {\u_{2,1}^{c}}^{T}Q \w_{1}, f_{12}= {\w_{2,1}^{c}}^{T} Q^{T} \u_{1}
\end{equation}
where 
\begin{eqnarray*}
\u_{1}:=\bmatrix{d_{1}\\d_{2}\\d_{3}}, 
\w_{1}:=\bmatrix{r_{1}\\r_{2}\\r_{3}},
\u_{1,1}^{c}=\bmatrix{d_{2}+d_{3}\\d_{1}+d_{3}\\d_{1}+d_{2}},
\u_{2,1}^{c}=\bmatrix{d_{2}d_{3}\\d_{1}d_{3}\\d_{1}d_{2}},
\w_{2,1}^{c}=\bmatrix{r_{2}r_{3}\\r_{1}r_{3}\\r_{1}r_{2}}.
\end{eqnarray*}
\textbf{Construction of Orthostochastic (Unistochastic) Matrices:}
From linear algebra it is known that if $\z_{i} \in \mathbf{col}(G_{i})$, column space of $G_{i}$ in equation (\ref{eqdiag}), the system is consistent. 
Here we have to study three cases separately.

\textbf{Diagonal matrices $D_{1}, D_{2}$ are simple}  (the eigenvalues of coefficient matrices $A_{i},i=1,2$ are all distinct): By Corollary \ref{Corollarydiag} $\diag(A_{i}, i=1,2$ are uniquely determined.
In order to compute a suitable orthostochastic matrix $Q$ we exploit the relations $\diag(A_{1})=Q^{T} \eig(A_{1})$ and $\diag(A_{2})=Q~ \eig(A_{2})$. The number of free parameters of $3\times 3$ doubly stochastic matrix is $4$. So, using one of the above two relations we can eliminate two diagonal entries by choosing two off-diagonal entries as parameters in the first $2 \times 2$ principal block sub-matrix of matrix $Q$. Similarly, using the other relation and imposing the condition that one orthostochastic matrix is the transpose of the other orthostochastic matrix, we can get rid of one more variable. 


Explicitly, 
\begin{eqnarray*}
 q_{11}&=\diag(A_{2})(0,0)-r_{3}-q_{12}(r_{2}-r_{3}))/(r_{1}-r_{3}) \\
 q_{21}&=\frac{(r_{1}-r_{3})(\diag(A_{1})(0,0)-d_{3})-(d_{1}-d_{3})(\diag(A_{2})(0,0)-r_{3}-q_{12}(r_{2}-r_{3})}{(r_{1}-r_{3})(d_{2}-d_{3})}\\
 q_{22}&=(\diag(A_{1})(1,0)-d_{3}-q_{12}(d_{1}-d_{3}))/(d_{2}-d_{3})
\end{eqnarray*}
where $\diag(A_{i})(0,0),\diag(A_{i})(1,0)$ denote the Ist and 2nd component of $\diag(A_{i})$. So, in the generic case of cubic bivariate we obtain a doubly stochastic matrix $Q=(q_{ij})$ in one parameter, say $q_{12}$. 

Note that there is a necessary and sufficient condition for a doubly stochastic matrix $Q$ of order $3$  to be an orthostochastic (unistochastic) matrix.  The condition is given by  \cite{Nakazato},\cite{Oleg}.
\begin{align} \label{eqortho}
(1-q_{11}-q_{12}-q_{21}-q_{22}+q_{11}q_{22}+q_{12}q_{21})^{2}=(\leq) 4q_{11}q_{22}q_{12}q_{21}
\end{align}
Using this necessary and sufficient condition we compute orthsostochastic matrices $Q$ and each of which provides us an orthogonal matrix $V$ such that $A_{2}=VD_{2}V^{T}$ and $f(\x)=\det(I+x_{1}D_{1}+x_{2}A_{2})$.

\noin \textbf{Degenerate Case}: Here we explain how to deal with degenerate cases which can be separated into two subcases. 

\noin At least one of two diagonal matrices satisfies the condition: \textbf{(All three diagonal entries are equal)}

Say wlog $D_{1} = \lambda I_{3}$, identity matrix of order $3$ and $\lambda$ is a non zero scalar, as it can't be a zero matrix. Observe that there are infinitely many (orthogonally equivalent) symmetric representations as 
\begin{equation*}
 f(\x)=\det(I+x_{1}\lambda I+x_{2}D_{2})=\det(I+x_{1}\lambda I+x_{2}VD_{2}V^{T})
\end{equation*}
for any orthogonal matrix $V$ of order $3$. 

\noin At least one of the two diagonal matrices satisfies the condition: \textbf{(Two of three eigenvalues are equal)} 
Wlog, say $D_{1}$ has two repeated eigenvalues. So, there are infinitely many ways to choose diagonal entries of coefficient matrix $A_{2}$. Any choice of diagonal entries of $A_{2}$ i.e., $\z_{2} \in$ col$(G_{2})$ in equation (\ref{eqdiag}) among infinitely many choices will work provided it is majorized by the vector $\u_{1}$, consisting of eigenvalues of $A_{2}$, see Example \ref{exdegenerate}.

\begin{remark}
 \rm{
 The existing methods in literature \cite{Helton}, \cite{Henrion}, \cite{Paprimulti} do not work in this case, but the method proposed in this paper enables us to compute a monic symmetric/Hermitian determinantal representation if it exists. Moreover, the method provides one representative candidate from each equivalence class of such a determinantal representation in generic case.
 }
 \end{remark}
Moreover we propose an algorithm which can efficiently compute such a determinantal polynomial for a cubic bivariate case.\\\\
\begin{algorithm}\label{algo2}
  \caption{Algorithm to Determine an MSDR of size $3$ by finding Orthostochastic Matrix}
\begin{algorithmic}
\small
\State Input: Cubic bivariate polynomial \begin{equation*}
 f(x_{1},x_{2})=
 f_{30}x_{1}^{3}+f_{03}x_{2}^{3}+f_{21}x_{1}^{2}x_{2}+f_{12}x_{1}x_{2}^{2}+f_{20}x_{1}^{2}+f_{02}x_{2}^{2}+f_{11}x_{1}x_{2}+f_{10}x_{1}+f_{01}x_{2}+1.
 \end{equation*}
\State Output: Orthostochastic matrix $Q=V \odot V$ such that
\begin{equation*}
f(\x)=\det(I+x_{1}A_{1}+x_{2}A_{2})=\det(I+x_{1}D_{1}+x_{2}VD_{2}V^{T})
\end{equation*}
\normalsize
\\\hrulefill
\begin{enumerate}
\item Determine the diagonal matrices $D_{1}, D_{2}$ by calculating the roots of univariate polynomials $f_{x_{1}}$ and $f_{x_{2}}$ respectively (Lemma \ref{lemmartdiag}).
\item Check that diagonal entries of $D_{1},D_{2}$ are real. If not, exit-no MSDR of size $3$ possible.
\item  Find the diagonal entries of $A_{1}$ and $A_{2}$  by solving two systems of linear equations defined in equation (\ref{eqdiag}).
\item Fix the descending order for the vectors $D_{1}=[\u_{1}]$ and $D_{2}=[\w_{1}]$. 
\item Check whether $\diag(A_{i}) \prec \eig(A_{i})$ for all $i=1,2$. If not, then MSDR of size $3$ is not possible-exit.
\item Find a doubly stochastic matrix $Q$ such that $Q^{T}\eig(A_{1})=\diag(A_{1}),$ and $Q\eig(A_{2}) =\diag(A_{2})$. Solution set is the intersection of Birkhoff Polytope $\B_{3}$ and a line segment. 
\item Find an orthostochastic matrix $Q$ such that $Q^{T}\eig(A_{1})=\diag(A_{1}),$ and $Q\eig(A_{2}) =\diag(A_{2})$ (use the necessary and sufficient condition for a doubly stochastic matrix of size $3$ to be an orthostochastic matrix). If no such orthostochastic matrix exists,  exit-no MSDR of size $3$ possible.
\item Find an orthogonal matrix $V=(v_{ij})$ such that $Q=(v_{ij}^{2})$ when orthostochastic matrix $Q$ is numerically known.
\item Construct $D_{1}$ and $A_{2}=VD_{2}V^{T}$.
\end{enumerate}
\end{algorithmic}
\end{algorithm} 
 
We explain the idea through examples.
\begin{example} \label{examcubicdet} \rm{
Consider the bivariate polynomial $f(x_{1},x_{2})=6x_{1}^{3}+36x_{1}^{2}x_{2}+66x_{1}x_{2}^{2}+36x_{2}^{3}+11x_{1}^{2}+42x_{1}x_{2}+36x_{2}^{2}+6x_{1}+11x_{2}+1.$
Vectors $\eig(A_{1}):=\u_{1}:=\bmatrix{3&2&1}^{T}, \eig(A_{2}):=\w_{1}=\bmatrix{6&3&2}^{T}, \diag(A_{1}):=\y_{1}:=\bmatrix{2.5&2&1.5}^{t} \prec \w_{1}$, and $\diag(A_{2}):=\y_{2}:=\bmatrix{4.5&4&2.5}^{T}$.
So, it's verified that majorization creteria are satisfied, i.e.,  $\diag(A_{i})\prec \eig(A_{i}), i=1,2$. Using the relations $\diag(A_{1})=Q^{T} \eig(A_{1})$ and $\diag(A_{2})=Q~ \eig(A_{2})$ and imposing the condition that one orthostochastic matrix is the transpose of the other one we obtain a line of doubly stochastic matrix $Q$ such that
\begin{equation*}
Q=\bmatrix{\frac{5-2u}{8} & u & \frac{3-6u}{8} \\ \frac{1+2u}{4} & 1-2u & \frac{6u-1}{4} \\
\frac{1-2u}{8} & u & \frac{7-6u}{8}}
\end{equation*} 
which satisfies the required conditions. It would be an orthostochastic if it satisfies the equation (\ref{eqortho}). 
So, we have a cubic equation $24u^{3}-6u+1=0$ which gives $u=-.5686,.3711,.1975$. Using the necessary and sufficient condition for an orthostochastic matrix $Q$ we find the values of $u=.37111, .19747$.

Note that if we know all the entries of the orthostochastic matrix $Q$ numerically, it is easy to find one possible orthogonal matrix $V$ such that $Q=V \odot V$ (see the DeterminantalRepresentations Package in \textit{Macaulay2}). For example, at $u=.37111$, one possible solution is
\tiny
\begin{align*}
Q \approx \bmatrix{.5322 & .3711 & .0967 \\ .4356 & .2578 & .3066 \\ .0322 & .3711 & .5967}, V \approx \bmatrix{
.7295 & .6092 & .3109 \\ -.6599 & .5077 & .5538\\ .1795 & -.6091 & .7724}, A_{2} \approx \bmatrix{4.5 & -1.6166 & 0.1527 \\-1.6166 & 4 & -0.7831 \\ 0.1527 & -0.7831 & 2.5}
\end{align*}
and at $u=.19747$, one possible solution is
\begin{equation*}
Q \approx \bmatrix{.5756 & .1975 & .2269\\.3487 & .6051 & .0462 \\.0756 & .1975 & .7269}, V \approx \bmatrix{-.7587 & .4444 &  .4763 \\.5905 & .7779 & .2150 \\.27501 & -.4444 & .8526}, A_{2} \approx \bmatrix{4.5 & -1.4465 & -1.0321 \\-1.4465 &  4 & 0.3040 \\-1.0321 & 0.3040 & 2.5}.
\end{equation*}
\normalsize
Note that if the coefficient of $x_{1}x_{2}^{2}$ is $63.9$,
$\y_{2}=\bmatrix{2.325&2.7&.975}^{T} \nprec \u_{1}$ and if the coefficient of $x_{1}x_{2}^{2}$ is $64$, $\y_{2} =\bmatrix{2.3333& 2.6667&1}^{T} \prec \u_{1}$, but no MSDR is possible. 

Also note that for a fixed vector coefficient $(f_{11},f_{21})$, the range of coefficient $f_{12}$ lies inside a closed interval. For example if we fix $(f_{11},f_{21})=(42,36)$, the coefficient $f_{12} \in [64.8, 66.8]$ are associated with a bivariate polynomial which has an MSDR, although the coefficient $f_{12} \in [64,68.9]$ satisfy the relation that diagonal entries of coefficient matrix $A_{12}$ is majorized by its eigenvalues. The Remark \ref{remarkindneccond} is verified.
}
\end{example}
\begin{example}\label{exdegenerate}\rm{(\textbf{Degenerate Case}) 
 Consider the polynomial
 \begin{equation*}
 f(\x)=162x_{1}^{3}-23x_{1}^{2}x_{2}+99x_{1}^{2}-8x_{1}x_{2}^{2}-10x_{1}x_{2}+18x_{1}+x_{2}^{3}-x_{2}^{2}-x_{2}+1
 \end{equation*}
Here $\u_{1}=\bmatrix{9\\6\\3}, \w_{1}=\bmatrix{1\\-1\\-1}$. So, $\diag(A_{2}):=Q \w_{1}=\bmatrix{-.7778 \\-.1111\\-.1111}$, uniquely determined by solving the system of equations defined in (\ref{eqdiag}).

On the other hand, $\diag(V^{T}D_{1}V):=Q^{T} \u_{1}=\bmatrix{5\\7\\6}$ (say) which is majorized by $\u_{1}$. Using the relation $Q \bmatrix{1\\-1\\-1}=\bmatrix{-.7778\\-.1111\\-.1111}$, we can write $Q=\bmatrix{.1111& u & .8889-u\\.4445 & v & .5555-v\\.4444 & 1-u-v& u+v-.4444}$. On the other hand using the relation $Q^{T} \bmatrix{9\\6\\3}=\bmatrix{5\\7\\6}$ we could derive a relation between $u$ and $v$ which is $6u+3v=4$. Therefore, $Q=\bmatrix{.1111& u & .8889-u\\.4445 & \frac{4-6u}{3} & 2u-.7777 \\ .4444 & u-1/3 & .8889-u}$. Now applying the necessary and sufficient condition for a doubly stochastic matrix to be an orthostochastic matrix, we get a quadratic equation as follows.
\begin{equation*}
1.8891 u^{2}-2.0743 u+.548785=0
\end{equation*}
Two solutions of this equation are $u=.4445,.6536$.
At $u=.4445$,
 \begin{equation*}
 Q=\bmatrix{.1111 & .4445 & .4444\\.4445 & .4444 & .1111\\.4445 & .1111 & .4444}, V=\bmatrix{-1/3 & 2/3 & 2/3 \\ 2/3 & -1/3 & 2/3 \\ 2/3 & 2/3 & -1/3}
 \end{equation*}
 and 
 \begin{equation*}
 f(\x) =\det(I+x_{1}\bmatrix{9 & 0 & 0 \\ 0 & 6 & 0 \\ 0 & 0 & 3}+x_{2}\bmatrix{-.7778 & -.4444 & -.44444 \\ -.4444 & -.11111 & .8889 \\ -.44444 & .8889 & -.1111}) \\
 \end{equation*}
 At $u=.6536$,
 \begin{align*}
 &Q=\bmatrix{.1111 & .6536 & .2353\\.4445 & .02613 & .52937\\.4444 & .324787 & .2353}, V=\bmatrix{.3333 & -.808455 & .485077 \\ -.6667 & .161648 & -.727578 \\ .66667 & .5699 & .48042} \mbox{and} \\
 &f(\x) =\det(I+x_{1}\bmatrix{9 & 0 & 0 \\ 0 & 6 & 0 \\ 0 & 0 & 3}+x_{2}\bmatrix{-.7778 & -.4444 & .44444 \\ -.4445 & -.1111 & .8889\\ .44444 & -.8889 & -.1111}) \\
 \end{align*}
 }
 \end{example}

\section{Computational Relaxation for Determinantal Polynomials} \label{secrelax}
We need to find a suitable orthostochastic matrix $Q$ (a unistochastic matrix) to determine an MSDR (MHDR) of size $d$ for a bivariate polynomial if it exists. More explicitly, while computing  an MSDR (MHDR) of size $d$ we need to apply a necessary and sufficient condition for which a doubly stochastic matrix of order $d$ would be an orthostochastic (unistochastic) matrix.

This problem is unresolved if the order of doubly-stochastic matrix is $\geq 4$. So  we propose a computational relaxation to the original problem
by finding a doubly stochastic matrix and its transpose such that they majorize both the doubly stochastic polytope in which diagonal entries are majorized by the eigenvalues of coefficient matrices instead of finding orthostochastic or unistochastic matrix. Moreover,
we conclude this subsection by providing a necessary and sufficient condition for the existence of such a doubly stochastic matrix.

After evaluating the values of vectors $\w_{1}$ and $Q\w_{1}$, we can get $d-1$ linear expressions in terms of entries of $Q:=(q_{ij})$, where $d$ is the size of doubly stochastic matrix $Q$. The number of free variables in a doubly stochastic matrix of size $d$ is $(d-1)^2$. So, we can eliminate $d-1$ free diagonal entries from $Q$ by parameterizing $(d-1)(d-2)$ off-diagonal entries of $M$ using the above mentioned $d-1$ linear expressions.

Note that there are two sets of $d-1$ monomials in which one of the two variables $x_{1},x_{2}$ must be of degree one, i.e., they are of the form $x_{1}^{\alpha}x_{2},$ or $x_{1}x_{2}^{\alpha}, \alpha \in \{1, \dots,d-1\}$ and monomial $x_{1}x_{2}$ is common in both the sets . So,  we can eliminate $d-2$ more free variables from the required doubly stochastic matrix $Q$  by using the values of $\u_{1}$ and $Q^{T}\u_{1}$. Therefore, the required doubly stochastic matrix $Q$ has $(d-1)(d-2)-(d-2)=(d-2)^{2}$ free off diagonal entries and it is a parameterized matrix in $(d-2)^{2}$ parameters in our context.

As each entries of $Q$ are linear in terms of these parameters and lies in the closed interval $[0,1]$ (due to the definition of doubly stochastic matrix), so we can specify feasible region for this system of linear multivariate inequalities in $(d-2)^{2}$ variables. Thus the problem turns into a problem of solving a system of linear multivariable inequalities. In Linear algebra Farkas Lemma (or theorem of the alternative) provides a certificate of emptyness for a polyhedral set $\{\x:A \x \leq \b \}$ for some matrix $A \in \R^{m \times n}$ and some vector $\b \in \R^{m}$. One can use the command LinearMultivariateSystem to solve a system of linear inequalities with respect to the given variables in Maple.

Solution of this system of linear multivariate system provides a tight region in which each doubly stochastic matrix and its transpose majorize both the required polytopes. In fact, if we relax the problem of determining orthostochastic (unistochastic) matrix  to a problem of determining a doubly stochastic matrix which satisfies the majorization criteria explained in  Proposition \ref{propnec}, it turns  into a problem of deciding whether a point lies inside a convex hull of  the finite set of specified points. Now by combining two necessary conditions mentioned in Proposition \ref{propnec} we provide a necessary and sufficient condition for the relaxation problem which is in fact a necessary condition for the original problem.

Consider the bivariate polynomial $f(x_{1},x_{2})$. Let $(f_{\alpha_{1},1},\dots,f_{1,\alpha_{2}})$ denotes the vector of coefficients of mixed monomials $x_{1}^{\alpha_{1}}x_{2},x_{1}x_{2}^{\alpha_{2}}, \alpha_{1},\alpha_{2} \in \{1,\dots,d-1\}$ of $f(x_{1},x_{2})$.
\begin{theorem} \label{themcommondoubly}  There exists a doubly stochastic matrix $Q$ such that the vector coefficient $(f_{\alpha_{1},1},\dots, f_{1,\alpha_{2}})$ of size $2d-2$ satisfies the scalar product representation mentioned in the Theorem \ref{themscalar} if and only if the vector coefficient $(f_{\alpha_{1},1},\dots,f_{1,\alpha_{2}})$ can be expressed as some convex combination of the following $d!$ points of size $2d-2$
\begin{align*}
\{{\u_{\alpha_{1},1}^{c}}^{T} P \w_{1}, {\w_{\alpha_{2},1}^{c}}^{T} (P)^{T} \u_{1},\alpha_{1},\alpha_{2}=1,\dots,d-1, P \
 \mbox{all permutation matrices of order} \ d\}
\end{align*}
where the vectors $\u_{1},\w_{1},\u_{\alpha_{1},1}^{c},\w_{\alpha_{2},1}^{c}$ are defined in equations (\ref{eqscalar1}) and (\ref{eqscalar2}).
\end{theorem}
\pf Suppose there exists such a doubly stochastic matrix $Q$. Then by the Theorem \ref{themscalar} the vector coefficient $(f_{\alpha_{1},1}, \dots, f_{1,\alpha_{2}})$ is of the form
\begin{equation} \label{eqrelaxation}
 (f_{\alpha_{1},1},\dots, f_{1,\alpha_{2}})=(\u_{\alpha_{1},1}^{c})^{T}Q \w_{1}, (\w_{\alpha_{2},1}^{c})^{T} Q^{T} \u_{1}), \ \forall \alpha_{1},\alpha_{2} \in \{1,\dots,d-1\}.
\end{equation}
Observe that the monomial $x_{1}x_{2}$ is common in  mixed monomials $x_{1}^{\alpha_{1}}x_{2},x_{1}x_{2}^{\alpha_{2}}, \alpha_{1},\alpha_{2} \in \{1,\dots,d-1\}$ and
\begin{equation*}
 f_{11}={\u_{1,1}^{c}}^{T} Q \w_{1}= {\w_{1,1}^{c}}^{T} Q^{T} \u_{1}.
\end{equation*}
On the other hand, it follows from the Theorem \ref{orthoconv2} that the range set $\{Q \w_{1}, Q \in \Omega_{d}\}$ is a convex set which is in fact a \textit{generalized permutohedron}. Using the property of linearity of second argument we have
\begin{equation*}
\{(1-\lambda) \langle \u,Q_{1}\w_{1} \rangle+ \lambda \langle \u,Q_{2}\w_{1} \rangle: Q_{1},Q_{2} \in \Omega_{n}\}=
\{\langle \u,Q\w_{1} \rangle, (1-\lambda)Q_{1}\w_{1}+\lambda Q_{2}=:Q \in \Omega_{n}\}
\end{equation*}
Thus, the set $\{\u^{T}Q \w_{1}=
\langle \u,\w_{1} \rangle_{Q}: Q \in \Omega_{d}\}$ is the convex hull of $\{\u^{T}P \w_{1}$ where $P$ is all possible permutation matrices of order $d$. As the Cartesian product of convex sets is a convex set, therefore, the set 
\begin{equation*}
 \{(\u_{\alpha_{1},1}^{c})^{T}Q \w_{1}, (\w_{\alpha_{2},1}^{c})^{T} Q^{T} \u_{1}\}, Q \in \Omega_{d} \ \forall \alpha_{1},\alpha_{2} \in \{1,\dots,d-1\} 
\end{equation*}
is a convex set. Moreover, it is the convex hull of $d!$ points of size $2d-2$ as follows.
\begin{align*}
\{{\u_{\alpha_{1},1}^{c}}^{T} P \w_{1}, {\w_{\alpha_{2},1}^{c}}^{T} (P)^{T} \u_{1},\alpha_{1},\alpha_{2}=1,\dots,d-1, P \
 \mbox{all permutation matrices of order} \ d\}
\end{align*}
Therefore, for a specific doubly stochastic matrix $Q$ which satisfies the equation (\ref{eqrelaxation}), the vector coefficient  $(f_{\alpha_{1},1},\dots,f_{1,\alpha_{2}})$ can be expressed as some convex combination of the specified points. 

Conversely, if there exists such a convex combination, that convex combination of the corresponding permutation matrices provides a doubly stochastic matrix which satisfies the conditions associated with the vector coefficient $(f_{\alpha_{1},1},\dots,f_{1,\alpha_{2}})$ by the Theorem \ref{themmaj}.

\qed

Note that each of these vector coefficients which is obtained by some convex combination of $d!$ specified points need not be
associated with a determinantal polynomial 
since that convex combination
of permutation matrices need not be an orthostochastic  or  a unistochastic matrix. Thus we have the following corollaries.
\begin{corollary} Consider the bivariate polynomial $f(x_{1},x_{2})$. There exists a doubly stochastic matrix $Q$ with a possible pair of coefficient matrices $(A_{1},A_{2})$ such that $Q^{T} \eig(A_{1})=\diag(A_{1})$ and $Q \eig(A_{2})=\diag(A_{2})$ if and only if the vector coefficient $(f_{\alpha_{1},1},\dots,f_{1,\alpha_{2}})$ 
of $f(x_{1},x_{2})$ can be expressed as some convex combination of the following $d!$ points
\begin{align*}
\{{\u_{\alpha_{1},1}^{c}}^{T} P \w_{1}, {\w_{\alpha_{2},1}^{c}}^{T} (P)^{T} \u_{1},\alpha_{1},\alpha_{2}=1,\dots,d-1, P \
 \mbox{all permutation matrices of order} \ d\}
\end{align*}
where the vectors $\u_{1},\w_{1},\u_{\alpha_{1},1}^{c},\w_{\alpha_{2},1}^{c}$ are defined in equations (\ref{eqscalar1}) and (\ref{eqscalar2}).
\end{corollary}

\begin{corollary}
If some convex combination of permutation matrices produce an orthostochastic (a unistochastic) matrix, then the same convex combination of the vector coefficient $(f_{11},\dots,f_{d-11},f_{11},f_{12},\dots,f_{1d-1})$ associated with corresponding permutation matrices provides a vector coefficient of mixed monomials of determinantal bivariate polynomial whose coefficient matrices belong to the same orbits.
\end{corollary}

So, expressing the vector coefficient of mixed  monomials $x_{1}^{\alpha_{1}}x_{2},x_{1}x_{2}^{\alpha_{2}}, \alpha_{1},\alpha_{2}=1,\dots,d-1$   as convex combination of $d!$ specified points is not a sufficient condition,
but this is a necessary condition for the existence of MSDR (MHDR) of size $d$ for bivariate polynomials and it shows a method to compute an MSDR (MHDR)  for higher degree ($\geq 4$) bivariate polynomials which need not be strictly RZ polynomials.

Eventually, we develop an algebraic combinatorial  method to determine an MSDR of size $d$ for a bivariate polynomial of degree in the next subsection. 
\subsection{Construction of Orthostochastic Matrices from Permutation Matrices}
In this subsection, we discuss a method to construct an orthostochastic matrix by using the properties of permutation matrices for quartic bivariate polynomial.
This idea leads us to get a heureistic method to compute determinantal representation for higher degree bivariate polynomials.

Note that the \textit{Grassmannian} $G(k,d)$, a smooth projective variety of dimension $k(d-k)$ embeds into $\P^{{d \choose k}-1}$. Each point of $G(k,d)$ corresponds to an $k$-dimensional linear subspace of a fixed $d$- dimensional vector space, although every point in $\P^{{d \choose k}-1}$ is not \textit{Grassmannian} point.

Say $m:=$$d \choose k$. First we talk about a special type of permutation matrices of order $m$
which are obtained as Hadamard product of $k$-th exterior power of permutation matrices of order $d$
with themselves, call them \textit{Grassmannian} permutation matrices. 

For example, there are $6!$ permutation matrices of order $6$ among which only $4!$ permutation matrices of order $6$ are \textit{Grassmannian} permutation matrices since they are obtained by taking Hadamard product of second exterior power of permutation matrices of order $4$ with themselves. 

This special type of permutation matrices, named as \textit{Grassmannian} permutation matrices play a crucial role to construct an orthostochastic matrix.

Using the properties of \textit{Grassmannian algebra}, we conclude that orthostochastic (unistochastic) matrices $Q^{\wedge^{k}}, k=1\dots,d,$ can be expressed as a
convex combination of \textit{Grassmannian} permutation matrices. For example, consider
\begin{equation*}
V=\bmatrix{1/\sqrt{3} & 0 & -\sqrt{2/3} & 0 \\0 & \sqrt{1/6} & 0 & -\sqrt{5/6}\\ \sqrt{2/3} & 0 & \sqrt{1/3} & 0 \\ 0 & \sqrt{5/6} & 0 & \sqrt{1/6}}.
\end{equation*}
Then
\begin{align*}
&Q=V \odot V=1/6 P_{1234}+2/3 P_{3412}+1/6P_{1432} \\
& \mbox{and} \\
&Q^{\wedge^{2}}=1/18P_{1234}^{\wedge^{2}} \odot P_{1234}^{\wedge^{2}} +5/9 P_{3412}^{\wedge^{2}} \odot P_{3412}^{\wedge^{2}}+5/18P_{1432}^{\wedge^{2}} \odot P_{1432}^{\wedge^{2}}+2/18P_{3214}^{\wedge^{2}} \odot P_{3214}^{\wedge^{2}}.
\end{align*}
where $\pi \in S_{n}$, the permutation (symmetric) group and $P_{\pi}$ is the corresponding permutation matrix.

The Bruhat graph of $S_{n}$ is the directed graph whose nodes are the elements of $S_{n}$
and whose edges are given by $x \rightarrow y$ means that $x \xrightarrow{t} y$ for some $t \in T$. Bruhat order is the partial order relation
on the set $S_{n}$ defined by the relation $x < y$ means that there exist adjacent transpositions $s_{i} = (i, i+1)$ such that
\begin{equation*}
x=x_{0} \rightarrow x_{1} \rightarrow \dots \rightarrow x_{k-1} \rightarrow x_{k}=y.
\end{equation*}
The Bruhat (strong) and right weak order graphs of symmetric group $S_{4}$ are shown in the following figure \cite{coxeter}.
\begin{figure}[h]
\begin{center}
\includegraphics[scale= 0.61]{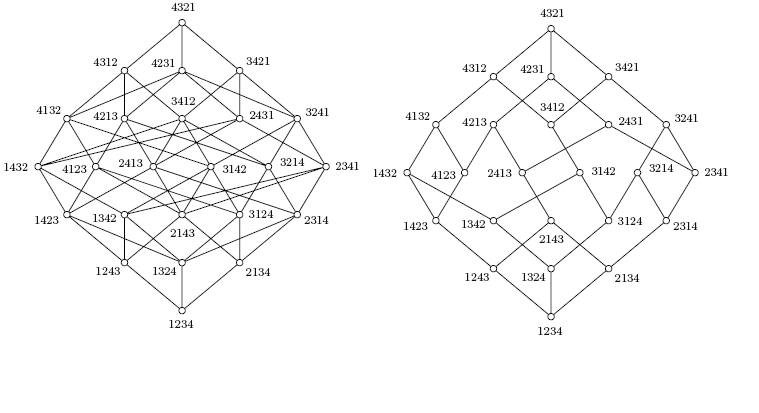}
\caption{Bruhat (Strong) and Weak order of $S_{4}$}\label{weakorder}
\end{center}
\end{figure}
As we can see in the Figure \ref{weakorder} , there are $4!=24$ permutations in the symmetric group $S_{4}$. Since the edges of the graph of (right) weak order of $S_{4}$ are constructed via an adjacent transposition, so we define a parameter $\alpha \in [0,1]$ for each point of the link of any two nodes $x,y$ (edges of) the graph $S_{4}$ such that $\alpha \leftrightarrow \alpha P_{x}+(1- \alpha)P_{y}$.

Thus, we have the following results about orthostochastic matrices of order $4$.
\begin{itemize}
\item Any convex combinations of any two edges of the Bruhat order graph of $S_{4}$ in Fig \ref{weakorder} correspond to orthostochastic matrices.
\item  Any convex combinations of any two adjacent nodes of the right weak order graph of $S_{4}$ in Fig \ref{weakorder} correspond to orthostochastic matrices.
\item The following collection of permutation matrices are such that any convex combination of four permutation matrices
from the set $B_{i} \in \mathcal{F}, i=1,\dots,18$  provide an orthostochastic matrix  in the Birkhoff polytope 
$\mathcal{B}_{4}$  due to its special structure.

$\mathcal{F}=\{B_{1}, \dots, B_{18}\}$ where 
\tiny
\begin{align*}
&B_{1}:=\{P_{2314},P_{3241},P_{1423},P_{4132}\}, B_{2}:=\{P_{2431},P_{1342},P_{3124},P_{4213}\}, 
B_{3}:=\{P_{2134}, P_{1243}, P_{4312}, P_{3421}\},\\
&B_{4}:=\{P_{3214}, P_{1432}, P_{2341}, P_{4123}\},
B_{5}:=\{P_{4231}, P_{1324}, P_{3142}, P_{2413}\}, B_{6}:=\{P_{1234},P_{2143},P_{3412},P_{4321},\}\\
&B_{7}:=\{P_{1234},P_{2134},P_{1243},P_{2143}\}, B_{8}:=\{P_{1324},P_{1342},P_{3142},P_{3124}\}, B_{9}:=\{P_{4231},P_{2431},P_{4213},P_{2413}\}, \\ &B_{10}:=\{P_{3412},P_{4321},P_{4312},P_{3412}\}, B_{11}:=\{P_{1432},P_{1423},P_{4123},P_{4132}\}, B_{12}:=\{P_{3214},P_{2314},P_{3241},P_{2341}\},\\
&B_{13}:=\{P_{2143},P_{3412},P_{3142},P_{2413}\}, B_{14}:=\{P_{3214},P_{3124},P_{4213},P_{4123}\}, B_{15}:=\{P_{2134},P_{2314},P_{4312},P_{4132}\},\\
 &B_{16}:=\{(P_{1432},P_{2431},P_{1342},P_{2341}), B_{17}:=\{P_{1234},P_{4231},P_{1324},P_{4321}\}, B_{18}:=\{P_{1243},P_{3241},P_{3412},P_{1423}\}
\end{align*}
\normalsize
\end{itemize}
Though it is complicated to get a figure of graph of $S_{n}$, but using the same arguments  we conclude that
any convex combinations of any two edges of the Bruhat order graph of $S_{n}$ correspond to orthostochastic matrices. Similarly, one can construct the regions which are the convex hulls of $d$ points  and  each point in Birkhoff polytope $\mathcal{B}_{n}$ is associated with an orthostochastic matrix.

However, we show how the relaxation method enables us to compute such a determinantal representation for lower degree cases in the following example. Also note that this method works better if we have repeated eigenvalues or in other words, polynomial is not strictly RZ polynomial.
\begin{example} \label{examquarbiv}
\rm{
Consider the quartic bivariate polynomial
\begin{align*}
f(x_{1},x_{2})&=24x_{1}^{4}+133.6609x_{1}^{3}x_{2}+50x_{1}^{3}+253.8824x_{1}^{2}x_{2}^{2}+
196.9412x_{1}^{2}x_{2}+35x_{1}^{2}\\&+190.4498x_{1}x_{2}^{3}+230.4498x_{1}x_{2}^{2}+87.6125x_{1}x_{2}+10x_{1}
+48x_{2}^{4}+80x_{2}^{3}\\&+48x_{2}^{2}+12x_{2}+1
\end{align*}
By Lemma \ref{lemmartdiag} the $\eig(D_{1})=\eig(A_{1}):=\u_{1}:=\bmatrix{d_{1}&d_{2}&d_{3}&d_{4}}^{T}=\bmatrix{4&3&2&1}^{T}$ and the $\eig(A_{2}):=\w_{1}:=\bmatrix{6&2&2&2}^{T}$.
Thus $\diag(A_{2}):=\y_{2}=\bmatrix{3.3840& 3.6749&   2.8858& 2.0554}^{T}$ by Corollary \ref{Corollarydiag}. By solving linear equations in entries of $Q$ of the form $Q  \w_{1}=\y_{2}$ we obtain the first column of matrix $Q$ as $\bmatrix{.346&.4187&.2215&.0138}^{T}$.
In order to determine second and third columns of matrix $Q$ we use the Theorem \ref{themcommondoubly}. Note that $\eig(A_{2})$ has three repeated eigenvalues, so by the Theorem \ref{themcommondoubly} the vector coefficients of monomials $x_{1}^{3}x_{2},x_{1}^{2}x_{2},x_{1}x_{2}^{3},x_{1}x_{2}^{2},x_{1}x_{2}$ can be expressed as convex combinations of $\frac{4!}{3!}=4$ specified points which are as follows
\begin{eqnarray*}
&\beta_{1}:=[132,196,192,232,88], \beta_{2}:=[124,184,176,216,84]\\
&\beta_{3}:=[148,216,208,248,92],\beta{4}:=[196,244,224,264,96]
\end{eqnarray*}
Note that $6$ permutation matrices are associated with each $\beta_{i},i=1,\dots,4$. Now our aim is to express the vector coefficient $[133.6609,196.9412,190.4498,230.4498,87.6125]$ as convex combination of these $4$ points such that the same convex combination of the corresponding permutation matrices provide an orthostochastic matrix.

Observe that one possible convex combination for the vector coefficient of given polynomial is
\begin{equation*}
[133.6609,196.9412,190.4498,230.4498,87.6125] 
=.4187\beta_{1}+.346 \beta_{2}+ .2215 \beta_{3}+.0138 \beta_{4}
\end{equation*}
Using the structure of $\mathcal{F}$ mentioned before we can choose 
\begin{equation*}
Q=.4187P_{2134}+.346P_{1243}+.2215P_{4312}+.0138P_{3421}\\
\end{equation*}
This convex combinations are found by solving  systems of linear equations.
This gives the orthostochastic matrix $Q=\bmatrix{0.3460 &0.4187& 0.0138& 0.2215\\   0.4187 &   0.3460 &  0.2215    &0.0138\\    0.2215 &   0.0138  &  0.4187  &  0.3460\\    0.0138 &   0.2215 &   0.3460 &   0.4187}$. Thus one possible orthogonal matrix $V$ and coefficient matrix $A_{2}$ are are as follows.
\begin{align*}
V \approx \bmatrix{.5882& .6471& .1175& .4706\\ .6471& -.5882& -.4706& .1174\\ .4706 &-.1174& .6471& -.5882\\ .1175& .4706 &-.5882& -.6471},A_{2} \approx \bmatrix{3.3840 & 1.5225 &1.1074 &  0.2764\\   1.5225  & 3.6748  & 1.2181&   0.3041\\ 1.1074 &  1.2181 &  2.886&  0.2211\\
   0.2764  & 0.3041 &  0.2211 &   2.0552}
\end{align*}
The pair of coefficient matrices $(D_{1},A_{2})$ satisfies the coefficient of monomial $x_{1}^{2}x_{2}^{2}$, so it provides a monic symmetric determinantal representation of the given polynomial.

}
\end{example}
\begin{remark}
\rm{
Observe that the quartic  bivariate polynomial in Example \ref{examquarbiv} is  not a  strictly RZ polynomial.
}
\end{remark}
\begin{remark}
\rm{
It is evident that if at least one coefficient matrix of a determinantal representation  has repeated eigenvalues, the given polynomial is not a strictly RZ polynomial. 
}
\end{remark}

There could be three possibilities
\begin{enumerate}
\item If the vector coefficient of the given polynomial cannot be expressed as convex combination of specified points, by the Theorem \ref{themcommondoubly} there is no such doubly stochastic matrix. This implies no such orthostochastic (unistochastic)  matrix exists, so conclude that MSDR (MHDR)  of size $d$ is not possible for the given bivariate polynomial.
\item Only one  doubly stochastic matrix exists. In this case that doubly stochastic matrix has to be orthostochastic (unistochastic)  matrix if an MSDR (MHDR)  exists for the given bivariate polynomial.
\item There are infinitely many doubly stochastic matrices. This does not ensure that there exists an orthostochastic (unistochastic)  matrix too, but it ensures the region of existence of orthostochastic (unistochastic)  matrix if it exists.
\end{enumerate}
 It is exemplified.
\begin{example}\rm{
Consider the bivariate polynomial $f(x_{1},x_{2})=6x_{1}^{3}+37.97x_{1}^{2}x_{2}+71.94x_{1}x_{2}^{2}+36x_{2}^{3}+11x_{1}^{2}+42.99x_{1}x_{2}+36x_{2}^{2}+6x_{1}+11x_{2}+1$.
As we have seen before $\u_{1}=\bmatrix{3&2&1}^{T}, \w_{1}=\bmatrix{6&3&2}^{T}$. There exists a doubly stochastic matrix $Q=\bmatrix{.5 & 0 & .5\\.5 &.01 & .49\\0 & .99 & .01}$ such that $Q\w_{1}=\bmatrix{4&4.01&2.99}^{T} \prec \w_{1}$ and $Q^{T}\u_{1}=\bmatrix{2.5&1.01&2.49}^{T} \prec \u_{1}$. There does not exist any  orthostochastic matrix along the line $u-2v+1=0$ in our method. Also note that the given polynomial is not a RZ polynomial since at $\x=(3,-1)$ its restricted univariate polynomial
has complex roots. So  the existence of doubly stochastic does not imply the existence of such an orthostochastic matrix.
}
\end{example}

\section{Range set of vector coefficients of mixed monomials} \label{secrange}
Let $S \R^{d \times d} (\H^{d \times d}(\C))$ be the space of all symmetric (Hermitian) matrices of order $d$. The Orbit of $A \in S \R^{d \times d}$ is
defined by $\mathcal{O}_{A}=\{ VAV^{T} : V \in O(d) \}$ and the orbit of $A \in H^{d \times d}(\C)$ is defined by
$\mathcal{O}_{A}=\{ UAU^{\ast} : U \in U(d) \}$.
By the Spectral theorem for symmetric (Hermitian) matrices it is clear that each of symmetric (Hermitian) matrices $A_{1}, A_{2}$
belongs to the unique orbit of a diagonal matrices $D_{1}$ and $D_{2}$,
denoted by $\mathcal{O}_{D_{1}}$, and $\mathcal{O}_{D_{2}}$ respectively.
Thus the vector space $S \R^{d \times d}(H^{d \times d}(\C))$ is a union of disjoint orbits of diagonal matrices.

Consider the class of bivariate polynomials $f(\x)$ having MSDR with coefficient matrices $A_{1}$ and $A_{2}$ which are obtained from the same orbits $\mathcal{O}_{D_{1}}$, and $\mathcal{O}_{D_{2}}$ respectively. Observe that  any two determinantal bivariate polynomials of the class $\{\det(I_{d}+x_{1}D_{1}+x_{2}VD_{2}V^{T})\}$ differ from each other by the vector coefficient of mixed monomials only.  In other words, the polynomials of this class share the same coefficients due to all monomials but mixed monomials. 
We say that they satisfy certain similarity pattern.

In this section, we would like to  classify all such determinantal bivariate polynomials which satisfy this similarity pattern. We provide a geometric structure of the range set of existence for vector coefficient of mixed monomials of bivariate polynomial of the this class. 

On the other hand, we are also interested to know if we replace a coefficient matrix $A_{j}$ by some arbitrary same type (symmetric /Hermitian)
matrix $\widehat{A}_j$ of same order such that the \textit{spectrums} of coefficient matrices $A_{j}$ and
$\widehat{A_{j}}$ are same; i.e., $\sigma_{A_{j}}=\sigma_{\widehat{A}_{j}}$, does there exist a relation between coefficients of $f(\x)$ and $\widehat{f}(\x)=\det(I+\sum_{j=1}^n \x_j\widehat{A}_j)?$ 

Let $S$ denotes the set of all coefficients of mixed monomials of bivariate polynomials which admit MSDR (MHDR)
of size $d$ with coefficient matrices $A_{1}, A_{2}$ belonging to the same orbits $\mathcal{O}_{D_{1}}$ and
$\mathcal{O}_{D_{2}}$ respectively.
So, by the Theorem \ref{themscalar} the set of ordered tuple
\begin{align} \label{eqrangeset}
S &=\{f_{\alpha_{1}\alpha_{2}}: \alpha_{1} \geq \alpha_{2}, 1 \leq \alpha_{1},  \alpha_{2} \leq d-1\} \times \{f_{\alpha_{1}\alpha_{2}}: \alpha_{1} \leq \alpha_{2}, 1 \leq \alpha_{1},  \alpha_{2} \leq d-1\} \nonumber \\
&=\{{\u_{\alpha_{1},\alpha_{2}}^{c}}^{T} Q^{\wedge^{\alpha_{2}}} \w_{\alpha_{2}},{\w_{\alpha_{2},\alpha_{1}}^{c}}^{T}
(Q^{\wedge^{\alpha_{1}}})^{T} \u_{\alpha_{1}}, 1 \leq \alpha_{1},\alpha_{2} \leq d-1\}
\end{align}
where $Q^{\wedge^{\alpha_{1}}}$ and $Q^{\wedge^{\alpha_{2}}}$ are orthostochastic (unistochastic) matrices.

We show that the set $S$ is not a convex set, but lies inside a convex hull of some finite known points. 
\begin{theorem} \label{convexhullthm2}
 Consider the class of bivariate polynomials having MSDR (MHDR) with coefficient matrices $A_{1},A_{2}$ coming
 from same orbits $\mathcal{O}_{D_{1}}, \mathcal{O}_{D_{2}}$ respectively. Then the set $S$ defined in equation (\ref{eqrangeset}) lies inside the convex hull of $H$ where
 \begin{equation*}
H=\{{\u_{\alpha_{1},\alpha_{2}}^{c}}^{T}((P)^{\wedge^{\alpha_{2}}})\w_{\alpha_{2}},  {\w_{\alpha_{2},\alpha_{1}}^{c}}^{T}((P)^{\wedge^{\alpha_{1}}})^{T}\u_{\alpha_{1}}, 1 \leq \alpha_{1},\alpha_{2} \leq d-1 \}
\end{equation*}
where $P$'s are all permutation matrices of size $d$.
\end{theorem}
\pf  If a bivariate polynomial of this class admits an MSDR (MHDR), by the Theorem \ref{themmaj} there exists a set of permutation matrices whose convex
combination give us the required orthostochastic (unistochastic) matrix $Q$ which
satisfies the vector coefficient of mixed monomials $x_{1}x_{2},\dots,x_{1}^{d-1}x_{2},x_{1}x_{2}^{2},\dots,x_{1}x_{2}^{d-1}$. Say
$Q=\alpha_{1}P_{1}+\dots +\alpha_{m}P_{m}$, where $\sum_{i}^{m}\alpha_{i}=1, \alpha_{i} \geq 0, 1 \leq m \leq d!$.
Note that 
$Q^{\wedge^{k}}
=(\alpha_{1}P_{1}+\dots+\alpha_{m}P_{m})^{\wedge^{k}}$ 
can be written as $\sum_{i=1}^{d!}\beta_{i}P_{i}^{\wedge^{k}}, \beta_{i} \geq 0$,
and $\sum_{i=1}^{d!}\beta_{i}=1$.
So, these orthostochastic (unistochastic) matrices
lie inside the convex hull of those permutation matrices which can obtained as some $k$-th exterior power of permutation matrices of size $d$. It follows from the Theorem \ref{orthoconv2} that the range set $\{\y \in \R^{n} | \y=Q \x, Q \in \Omega_{n}\}$ is
a convex set. Thus, the set $\{\u^{T}Q \w_{1}=
\langle \u,\w_{1} \rangle_{Q}: Q \in \Omega_{n}\}$ is a convex set for any vectors $\u,\w_{1} \in \R^{n}$
such that $\u=\diag(D)$. But the set 
$\{\u^{T}Q\w_{1}, \w^{T}Q^{T} \u_{1}: Q \ \mbox{is an orthostochastic (unistochastic) matrix}\}$ is not a convex set. Thus the set $S$ is not a convex set. In fact, by Birkhoff Von Neumann theorem it is known that  
the set of all $d \times d$ doubly stochastic matrices $\Omega_{d}$ is the convex hull of $d \times d$ permutation matrices. The set of all $d\times d$ orthostochastic (unistochastic) matrices lies inside $\Omega_{n}$, but forms a nonconvex set. Using the same line of thoughts we conclude that $S$ lies 
inside the convex hull of $H$. Note that $H$ has $d!$ points.  \hfill{$\square$}

We show that the set $S$ attains its maximum and minimum at some specified points using the result of rearrangement inequality.

\noin Rearrangement inequality states that
\begin{align*}
&x_{n}y_{1}+\cdots +x_{1}y_{n}\leq x_{\sigma (1)}y_{1}+\cdots +x_{\sigma (n)}y_{n}\leq x_{1}y_{1}+\cdots +x_{n}y_{n} \\
&\mbox{for every choice of real numbers} \
x_{1}\leq \cdots \leq x_{n}\quad {\text{and}}\quad y_{1}\leq \cdots \leq y_{n} \\
&\mbox{and every permutation} \
x_{\sigma (1)},\dots ,x_{\sigma (n)}.
\end{align*}

\begin{proposition}
 Consider the class of bivariate polynomials having MSDR (MHDR) with coefficient matrices $A_{1},A_{2}$ coming
 from same orbits $\mathcal{O}_{D_{1}}, \mathcal{O}_{D_{2}}$ respectively. The set $S$ defined in equation (\ref{eqrangeset}) attains its minimum when the vectors $\u_{1}:=\diag(D_{1})$ and $\w_{1}:=\diag(D_{2})$  are in descending order and attains its maximum  when one of the  vectors $\u_{1}$ and $\w_{1}$ is  in descending order and the other one is in ascending order.
\end{proposition}
\pf  If both the vectors $\u_{1},\w_{1}$ are in descending order, then the vectors $\u_{k}, \w_{k}$ are in descending and $\u_{k,1}^{c}$ and $\w_{k,1}^{c}$ are in ascending order. From the Theorem \ref{themscalar}, it is clear that if the defining matrix is the identity permutation matrix associated with permutation $[1 \ 2 \  \dots \ n]$, then the coefficients of mixed monomials of a bivariate polynomial are scalar product of two vectors, one of which is descending and other one is ascending. Therefore, by the result of rearrangement inequality, it provides the minimum value of the set $S$. On the other hand, if the defining matrix is the permutation matrix associated with the permutation $[n \dots  2 \ 1]$, then it provides the maximum value of the set $S$ since the coefficients are obtained as scalar product of two descending vectors. \hfill{$\square$}

\begin{remark}
\rm{
All these results hold for monic Hermitian determinantal representation too.
}
\end{remark}
The result of the Theorem \ref{convexhullthm2} reflects the reason behind the Remark \ref{remarkindneccond}


\bibliographystyle{alpha}
\bibliography{refqplmi}

\begin{thebibliography}{GKVVW14}

\bibitem[BB06]{coxeter}
Anders Bjorner and Francesco Brenti.
\newblock {\em Combinatorics of Coxeter groups}, volume 231.
\newblock Springer Science \& Business Media, 2006.

\bibitem[BGFB]{Boydlmi}
S~Boyd, LE~Ghaoui, E~Feron, and V~Balakrishnan.
\newblock Linear matrix inequality in system and control theory. 1994.
\newblock {\em SIAM, Philadelphia, PA}, pages 7--35.

\bibitem[Br{\"a}10]{Brahyper}
Petter Br{\"a}nd{\'e}n.
\newblock Notes on hyperbolicity cones, 2010.

\bibitem[Bru06]{Brualdi}
Richard~A Brualdi.
\newblock {\em Combinatorial matrix classes}, volume~13.
\newblock Cambridge University Press, 2006.

\bibitem[CD08]{Oleg}
Oleg Chterental and Dragomir~{\v{Z}} Dokovi{\'c}.
\newblock On orthostochastic, unistochastic and qustochastic matrices.
\newblock {\em Linear Algebra and Its Applications}, 428(4):1178--1201, 2008.

\bibitem[Dey]{Paprimulti}
Papri Dey.
\newblock Definite determinantal representations of multivariate polynomials.
\newblock {\em https://arxiv.org/abs/1708.09557}.

\bibitem[Dix02]{Dixon}
Alfred~Cardew Dixon.
\newblock Note on the reduction of a ternary quantic to a symmetrical
  determinant.
\newblock In {\em Proc. Cambridge Philos. Soc}, volume~5, pages 350--351, 1902.

\bibitem[DP18]{paprirz}
Papri Dey and Daniel Plaumann.
\newblock Testing hyperbolicity of real polynomials.
\newblock {\em arXiv preprint arXiv:1810.04055}, 2018.

\bibitem[GKVVW14]{Vinnikov2}
Anatolii Grinshpan, Dmitry~S Kaliuzhnyi-Verbovetskyi, Victor Vinnikov, and
  Hugo~J Woerdeman.
\newblock Stable and real-zero polynomials in two variables.
\newblock {\em Multidimensional Systems and Signal Processing}, 27(1):1--26,
  2014.

\bibitem[Hen10]{Henrion}
Didier Henrion.
\newblock Detecting rigid convexity of bivariate polynomials.
\newblock {\em Linear Algebra and its Applications}, 432:1218--1233, 2010.

\bibitem[Hor54]{Horn}
Alfred Horn.
\newblock Doubly stochastic matrices and the diagonal of a rotation matrix.
\newblock {\em American Journal of Mathematics}, 76:620--630, 1954.

\bibitem[HV07]{Helton}
J.W. Helton and Vinnikov.
\newblock Linear matrix inequality representation of sets.
\newblock {\em Communications on Pure and Applied Mathematics}, 60:654--674,
  2007.

\bibitem[JT18]{Thorstenhyp}
Thorsten J{\"o}rgens and Thorsten Theobald.
\newblock Hyperbolicity cones and imaginary projections.
\newblock {\em Proceedings of the American Mathematical Society}, 2018.

\bibitem[KPV15]{Vinzant}
Mario Kummer, Daniel Plaumann, and Cynthia Vinzant.
\newblock Hyperbolic polynomials, interlacers, and sums of squares.
\newblock {\em Mathematical Programming}, 153(1):223--245, 2015.

\bibitem[LP17]{Plaumann}
Anton Leykin and Daniel Plaumann.
\newblock Determinantal representations of hyperbolic curves via polynomial
  homotopy continuation.
\newblock {\em Mathematics of Computation}, 86(308):2877--2888, 2017.

\bibitem[LPR05]{Pablo3}
A.~S. Lewis, P.~A. Parrilo, and M.~V. Ramana.
\newblock The lax conjecture is true.
\newblock {\em Proceedings of The American Mathematical Society},
  133:2495--2499, 2005.

\bibitem[MOA11]{Marshalinequality}
Albert~W. Marshall, Ingram Olkin, and Barry~C. Arnold.
\newblock {\em Inequalities: Theory of Majorization and Its Applications}.
\newblock Springer, 2011.

\bibitem[Nak96]{Nakazato}
Hiroshi Nakazato.
\newblock Set of 3 x 3 orthostochastic matrices.
\newblock {\em Nihonkai Math.J}, 7:83--100, 1996.

\bibitem[NN94]{Nesterov}
Y.~Nesterov and A.~Nemirovski.
\newblock {\em Interior-point polynomial algorithms in convex programming, vol.
  13 of SIAM Studies in Applied Mathematics}.
\newblock Society for Industrial and Applied Mathematics (SIAM), Philadelphia,
  PA, 1994.

\bibitem[Nt12]{Tim}
Tim Netzer and Andreas thom.
\newblock Polynomials with and without determinantal representations.
\newblock {\em Linear Algebra and its Applications}, 437:1579--1595, 2012.

\bibitem[Nui69]{Nuij}
W.~Nuij.
\newblock A note on hyperbolic polynomials.
\newblock {\em Math. Scand.}, 23:69--72, 1969.

\bibitem[Pav17]{Vincentext}
Vincent Pavan.
\newblock {\em Exterior Algebras}.
\newblock Elsevier, 2017.

\bibitem[PSV12]{Sturmfelsbivariate}
Daniel Plaumann, Bernd Sturmfels, and Cynthia Vinzant.
\newblock Computing linear matrix representations of helton-vinnikov curves.
\newblock In {\em Mathematical methods in systems, optimization, and control},
  pages 259--277. Springer, 2012.

\bibitem[Ram95]{Ramana1}
Motakuri Ramana.
\newblock Some geometric results in semidefinite programming.
\newblock {\em Journal of Global Optimization}, 7:33--50, 1995.

\bibitem[RRS]{Nikhilhyp}
Prasad Raghavendra, Nick Ryder, and Nikhil Srivastava.
\newblock Real stability testing.

\bibitem[TM01]{Taylor}
Michael~D. Taylor and Piotr Mikusinski.
\newblock {\em An Introduction to Multivariable Analysis from Vector to
  Manifold}.
\newblock Springer, 2001.

\bibitem[Win10]{Sergeiexterior}
Sergei Winitzki.
\newblock {\em Linear algebra via exterior products}.
\newblock Sergei Winitzki, 2010.

\end{thebibliography}
\end{document}